\documentclass[review]{elsarticle}
\makeatletter
\def\ps@pprintTitle{%
	\let\@oddhead\@empty
	\let\@evenhead\@empty
	\def\@oddfoot{}%
	\let\@evenfoot\@oddfoot}
\makeatother
\usepackage[margin=1in]{geometry}
\usepackage{float}
\usepackage[hidelinks]{hyperref}
\usepackage{dsfont}
\usepackage{bbold}
\usepackage{soul}
\usepackage{color}
\usepackage{mathrsfs}
\usepackage{hyperref}
\usepackage{url}
\usepackage{amsmath}
\usepackage{amssymb}
\usepackage{enumitem}
\usepackage{algorithm} 
\usepackage{algpseudocode}
\usepackage{mathtools}
\usepackage{comment}
\usepackage{url} 
\usepackage{soul}
\usepackage{subcaption}
\usepackage{stmaryrd}
\usepackage{amsthm}

\usepackage{pgfplots}
\usetikzlibrary{spy}
\newlength{\figSize}
\usepackage{graphicx} % \scalebox
\usepackage{environ}
\NewEnviron{myequation}{%
\begin{equation}
\scalebox{1.4}{$\BODY$}
\end{equation}
}
\begin{document}
	
	\begin{frontmatter}
		
		\title{Isoparametric Tangled Finite Element Method for Nonlinear Elasticity}
		
	
		\author[]{Bhagyashree Prabhune}
		\ead{bprabhune@wisc.edu}		
		\author[]{Krishnan Suresh\corref{mycorrespondingauthor}}
		\cortext[mycorrespondingauthor]{Corresponding author}
		\ead{ksuresh@wisc.edu}
		\address{Department of Mechanical Engineering, University of Wisconsin-Madison, WI, USA}
		
		

\begin{abstract}
%Elements with negative Jacobian i.e. tangled or concave elements are commonly encourerd during mesh generation.   
An important requirement in the standard finite element method (FEM) is that all elements in the underlying mesh must be tangle-free i.e., the Jacobian must be positive throughout each element. 
To relax this requirement, an isoparametric tangled finite element method (i-TFEM) was recently proposed for linear elasticity problems.  It was demonstrated that i-TFEM leads to optimal convergence even for severely tangled meshes.  
 
In this paper, i-TFEM is generalized to nonlinear elasticity. Specifically, a variational formulation is proposed that leads to local modification in the tangent stiffness matrix associated with tangled elements, and an additional piece-wise compatibility constraint. i-TFEM reduces to standard FEM for tangle-free meshes. The effectiveness and convergence characteristics  of i-TFEM are demonstrated through a series of numerical experiments, involving both compressible and in-compressible problems. 

\end{abstract}
		
		
	\end{frontmatter}
	
	%\linenumbers
	
\section{Introduction}	

The finite element method (FEM) is extensively used for solving a wide variety of problems in nonlinear solid mechanics. A crucial step in FEM is meshing, where a finite element mesh is generated to represent the physical domain. An important requirement of a valid FEM mesh is that it should not contain tangled elements \cite{zienkiewicz2005finite, frey2007mesh, cook2007concepts, lo2014finite}. 
In other words, every element must be fully invertible. However, generating high-quality, tangle-free meshes remains a challenge, despite advances in mesh generation methods \cite{blacker2001automated, pietroni2022hex}. State-of-the-art mesh generation methods \cite{mandad2022intrinsic, fang2016all, livesu2013polycut, gregson2011all, jiang2013frame, li2012all, nieser2011cubecover, huang2011boundary} can still lead to a tangled mesh.  As a result, many untangling techniques have been developed \cite{livesu2015practical, xu2018hexahedral, knupp2001hexahedral, ruiz2015simultaneous,  huang2022untangling}, but they are not always reliable. Several instances have been reported where no tangle-free solution could be found without  altering the boundary \cite{livesu2015practical, akram2021embedded, xu2018hexahedral}.
 
This paper is focused on handling inverted (i.e., concave) bilinear quadrilateral (Q4) elements, a widely used element in FEM, within the context of nonlinear solid mechanics. Various  non-traditional finite element formulations have been proposed to handle concave element shapes; for instance, smoothed finite element \cite{liu2007smoothed}, polygonal finite elements \cite{sukumar2004conforming}, virtual elements  \cite{beirao2013basic}, unsymmetrical FEM \cite{cen2015unsymmetric}. These have been applied to solve problems with finite deformations, for instance, PolyFEM \cite{chi2015polygonal, bishop2020polyhedral}, VEM \cite{wriggers2017efficient, van2020virtual, chi2017some}, unsymmetrical FEM \cite{li2020hyperelastic}. However, these methods often require significant modifications to FEM (for example, the use of non-standard shape functions) and/or  do not reduce to standard FEM for non-tangled (regular) meshes.

%Various non-traditional finite element methods have been proposed to handle tangled elements and have been applied to solve the problems in nonlinear elasticity. Methods such as the smoothed finite element method \cite{liu2007smoothed}, polyonal finite elements, and virtual element method \cite{beirao2013basic} could be employed to handle tangled (concave) elements and have been applied for finite deformation problems, for instance, see \cite{wriggers2017efficient, van2020virtual, de2019serendipity} for VEM, \cite{chi2015polygonal} for polyFEM, \cite{} for SFEM. Unsymmetric finite elements \cite{cen2015unsymmetric} have also been proposed to handle tangled elements; their application to finite deformations is reported in \cite{}.  However, as the name suggests, they result in unsymmetric stiffness matrices. Moreover, these methods do not reduce to standard FEM for non-tangled (regular) meshes.

Recently, an isoparametric tangled finite element method (i-TFEM) was proposed as an extension to the standard FEM \cite{prabhune2023computationally} to handle tangled elements for linear problems. The proposed method uses the same shape functions as FEM, reduces to the standard FEM for non-tangled meshes, and maintains the symmetry of the stiffness matrix.   In i-TFEM, tangled elements are handled by: (1) locally modifying the elemental stiffness matrices associated with the tangled elements,  and (b) enforcing  piecewise compatibility constraints at  re-entrant nodes.  It was demonstrated in the context of tangled quadrilateral (Q4) and 8-node (H-8) hexahedral meshes for linear elasticity and Poisson problems \cite{prabhune2023computationally, prabhune2022towards, prabhune2022tangled}. 

In this paper, i-TFEM is extended to handle finite deformation problems, with  hyperelasticity, over tangled 4-node quadrilateral (Q4) meshes.  Towards this end, a mixed formulation  for i-TFEM is introduced, and the required modifications to the tangent matrices over the concave elements are formally derived. It is demonstrated that the results obtained over tangled meshes using i-TFEM have the convergence rate comparable to that of the standard FEM with regular meshes. Several examples with compressibility and near-incompressibility are presented to demonstrate the robustness of i-TFEM. 

The remainder of this paper is organized as follows. The standard finite element  formulation for nonlinear elasticity is reviewed in Section 2; the challenges associated with tangled meshes are highlighted. Section 3 describes the proposed i-TFEM method using total Lagrangian formulation. This is followed by numerical experiments in Section 4, and conclusions in Section 5.

%%%%%%%%%%%%%%%%%%%%%%%%
\section{Nonlinear elasticity}
Consider a body occupying a domain $\Omega \in \mathbb{R}^2$ subject to a body force $\boldsymbol{b}$, traction $\boldsymbol{T}$ over the boundary $\partial \Omega ^T$, and Dirichlet boundary conditions $\boldsymbol{u} = \boldsymbol{u_d}$ over the boundary $\partial \Omega ^d$; the  material is assumed to be hyper-elastic undergoing a finite deformation $\boldsymbol{u}$. The domain is divided into $M$ elements $E_j$, identified by the set $I = \left\lbrace 1,\dots M\right\rbrace $. We employ the total Lagrangian formulation \cite{reddy2014introduction} in this paper, where the potential energy can be written as:
\begin{equation}
	\Pi(\boldsymbol{u}) = \sum_{j \in I} \int\limits_{E_j}{\Psi\left(\boldsymbol{F}\left(\boldsymbol{u}_j \right)  \right) dV}																						 - 																						\sum_{j \in I} \int	\limits_{E_j}	{\boldsymbol{u}_j \cdot \boldsymbol{b} dV} 																									- 																					\sum_{j \in I} \int\limits_{\partial E_j^T}{\boldsymbol{u}_j\cdot\boldsymbol{T}dS} 
	\label{Eq_Pi}
\end{equation}
where $\boldsymbol{F}$ is the deformation gradient, $\Psi$  is the strain energy density. Further, using the standard (Bubnov-) Galerkin variational formulation, one arrives at the residual equation \cite{reddy2014introduction}:
\begin{equation}
	\boldsymbol{R}(\boldsymbol{\hat{u}})  = 0
	\label{Eq_R_FEM}
\end{equation}
This is typically solved iteratively via the  Newton-Raphson algorithm \cite{reddy2014introduction}:
\begin{equation}
	\boldsymbol{K} (\boldsymbol{\hat{u}}^n) \Delta \boldsymbol{\hat{u}}^{n+1} = - \boldsymbol{R}(\boldsymbol{\hat{u}}^n).
	\label{Eq_FEMFiniteStrain}
\end{equation}
where $\boldsymbol{K}$ is the tangent  matrix and $\Delta\boldsymbol{\hat{u}}^n$ is the incremental displacement vector at $n^{th}$ Newton iteration. When the mesh is of high-quality and not tangled, one obtains accurate solutions to such problems. 

\begin{figure}[H]
	\begin{subfigure}[c]{.5\textwidth}				
		\centering\includegraphics[width=0.65\linewidth]{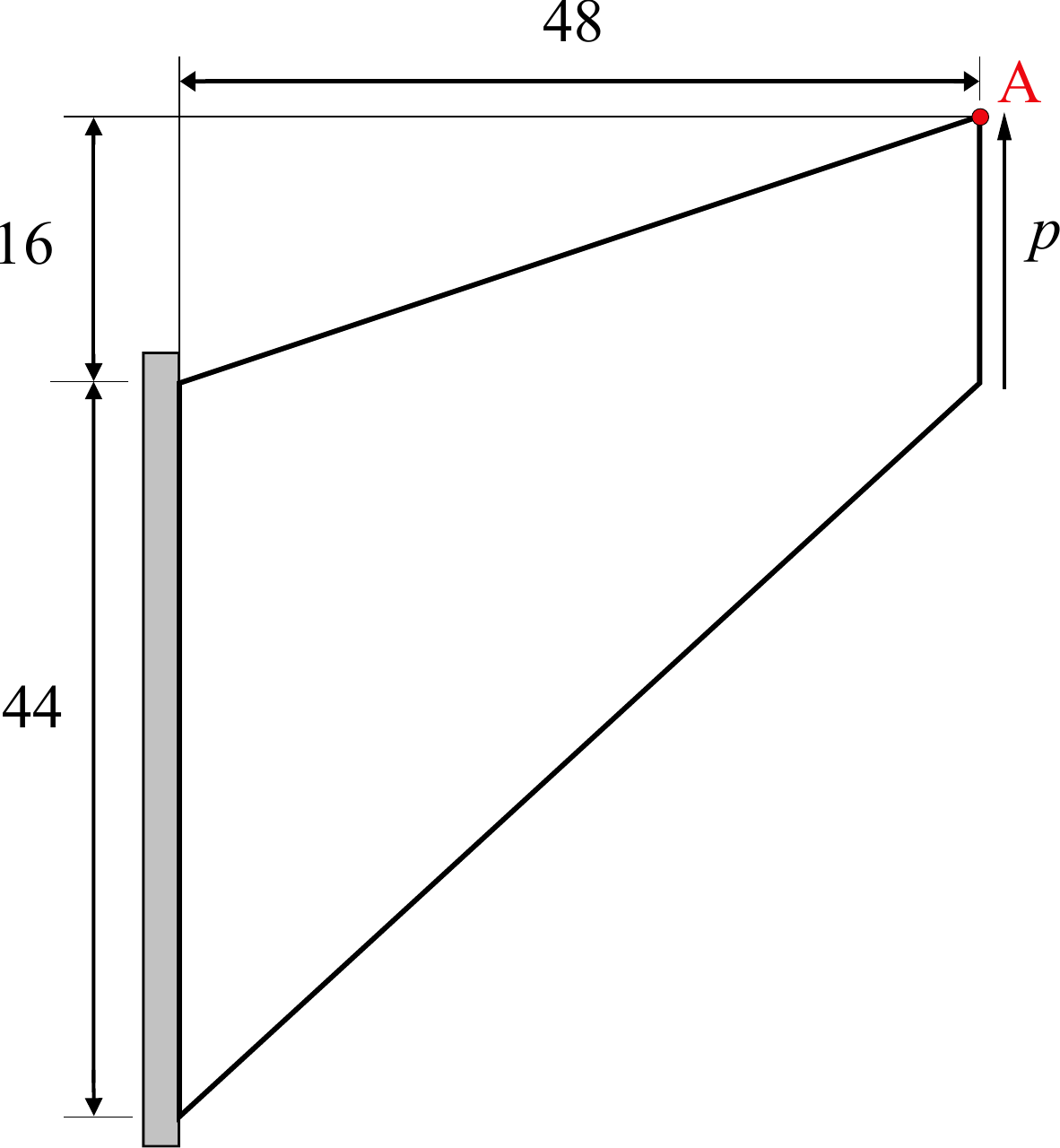}
		\caption{}
		\label{Fig_cooksGeom}
	\end{subfigure}
	\begin{subfigure}[c]{.5\textwidth}				
		\centering\includegraphics[width=0.55\linewidth]{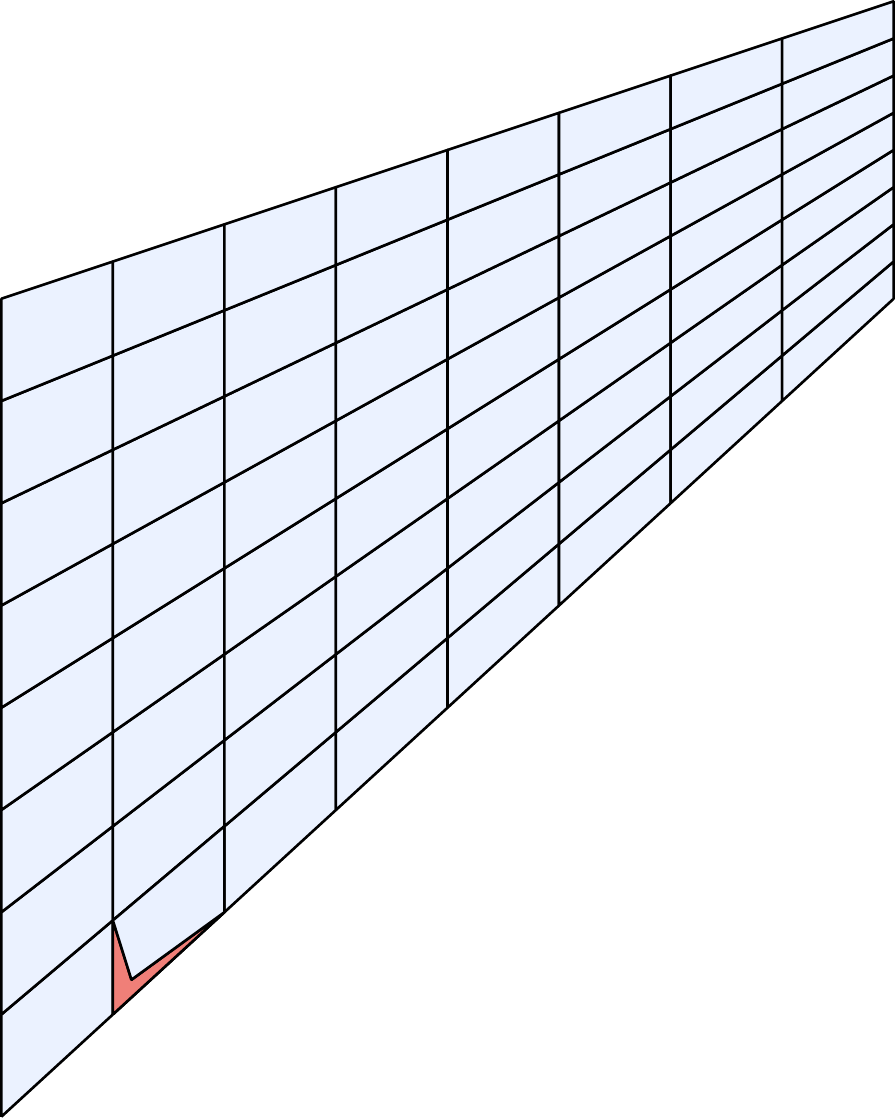}
		\caption{}
		\label{Fig_cooksMesh1tangle}
	\end{subfigure}	

	\caption{(a) Cook's membrane problem. (b) Tangled mesh with one concave element. }
	\label{}	
\end{figure}

\begin{figure}[H]
	\begin{subfigure}[c]{.5\textwidth}				
		\centering\includegraphics[width=0.5\linewidth]{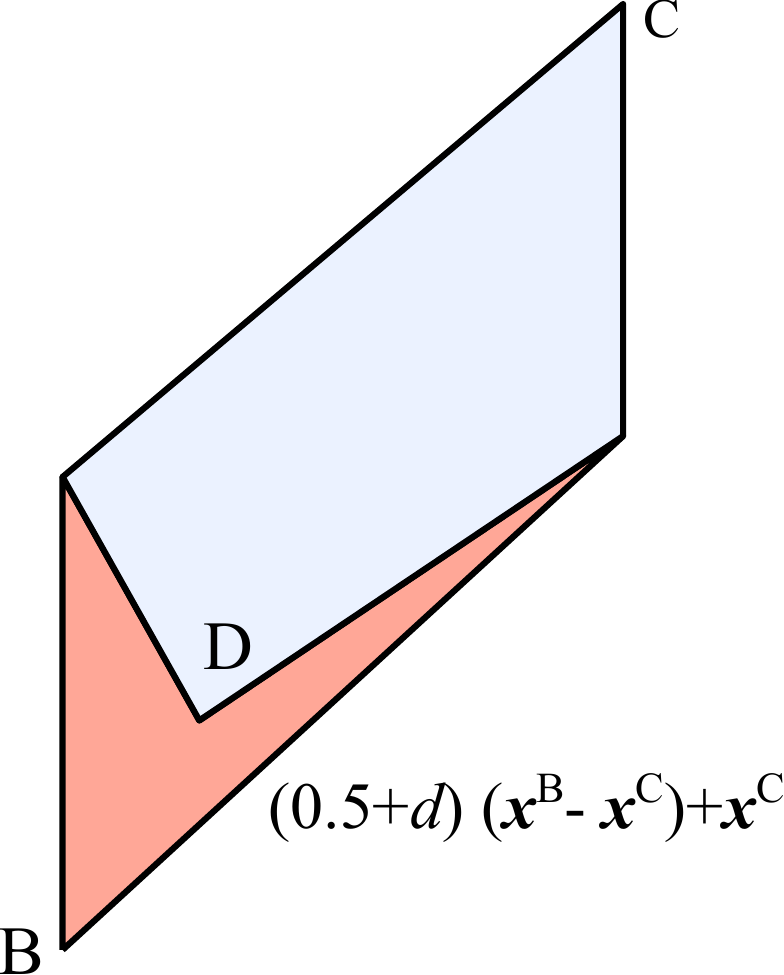}
		\caption{}
		\label{Fig_cook1Elem}
	\end{subfigure}
	\begin{subfigure}[c]{.5\textwidth}		
		\centering\includegraphics[width=0.8\linewidth]{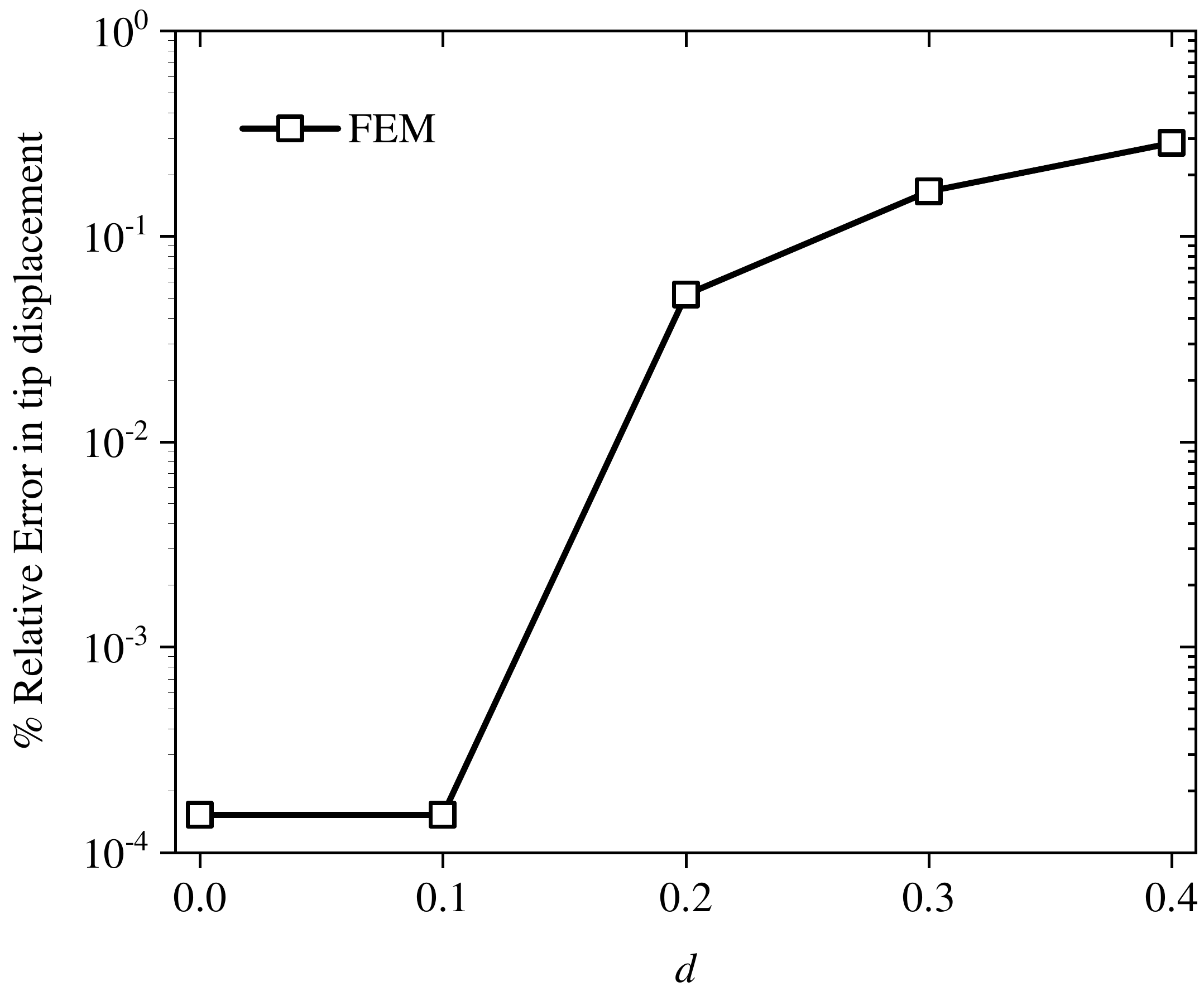}
		\caption{}
		\label{Fig_cooks1ElemResultsFEM}
	\end{subfigure}	
	\caption{(a) Zoomed-in view of the tangled element. (b) Relative error in tip displacement versus $d$ for FEM.}
	\label{}	
\end{figure}

However, as is well known, when the mesh is tangled, i.e.,  if the mesh contains inverted elements, the solution becomes erroneous. To illustrate, consider  Cook's membrane problem \cite{cook1974improved} illustrated in Fig.~\ref{Fig_cooksGeom}. The left edge of the tapered cantilever is fixed while a uniformly distributed load $p = 5$ is applied on the right edge.  We pose a   geometrically nonlinear plane-strain problem with Lamé parameters $\mu = 50$ and $\lambda = 100$.  Fig.~\ref{Fig_cooksMesh1tangle} illustrates a quadrilateral mesh with one concave element that we use for this experiment.
We vary the extent of tangling by moving the re-entrant vertex D along the diagonal BC as shown in Fig.~\ref{Fig_cook1Elem}.  When the parameter $d=0$, the point D lies half-way between B and C, and when $0< d<0.5$, the point D moves towards B, i.e., the element gets tangled. The large-deformation problem is solved using the normal procedure as described above, with 10 load steps. The tip displacement is compared against the expected value (using a high quality non-tangled mesh).  When $d> 0.1$, a  negative $|\boldsymbol{J}|$ value is encountered at one or more Gauss points, and Fig.~\ref{Fig_cooks1ElemResultsFEM} illustrates the resulting erroneous solution.

\section {Isoparametric TFEM for nonlinear elasticity}
\label{Section3_TFEM}	

The objective of this paper is to propose an iso-parametric tangled finite element method (i-TFEM), as a simple extension to classic FEM, for solving large deformation problems over tangled meshes. As a background, we briefly review the critical i-TFEM concepts proposed in  \cite{prabhune2022tangled} for linear problems.

\subsection{Isoparametric TFEM}

\begin{figure}[H]
	\begin{subfigure}[c]{.20\textwidth}				
		\centering\includegraphics[width=1\linewidth]{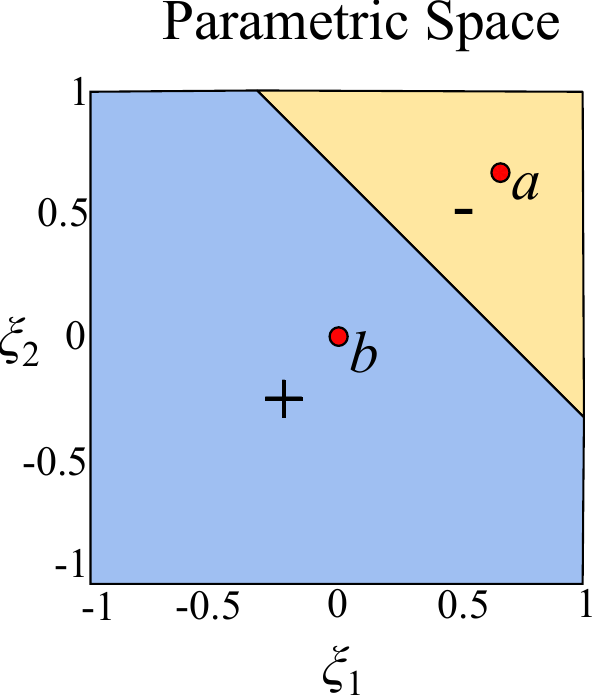}
		\caption{}
		\label{Fig_ConcaveQualMap_a}
	\end{subfigure}	
	\begin{subfigure}[c]{.20\textwidth}				
		\centering\includegraphics[width=1.1\linewidth]{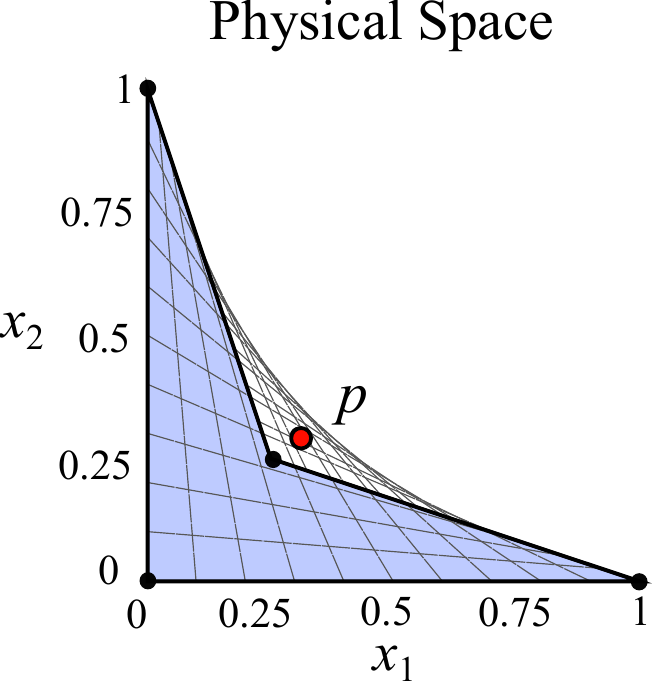}
		\caption{}
		\label{Fig_ConcaveQualMap_b}
	\end{subfigure}	
	\begin{subfigure}[c]{.6\textwidth}		
		\centering\includegraphics[width=0.9\linewidth]{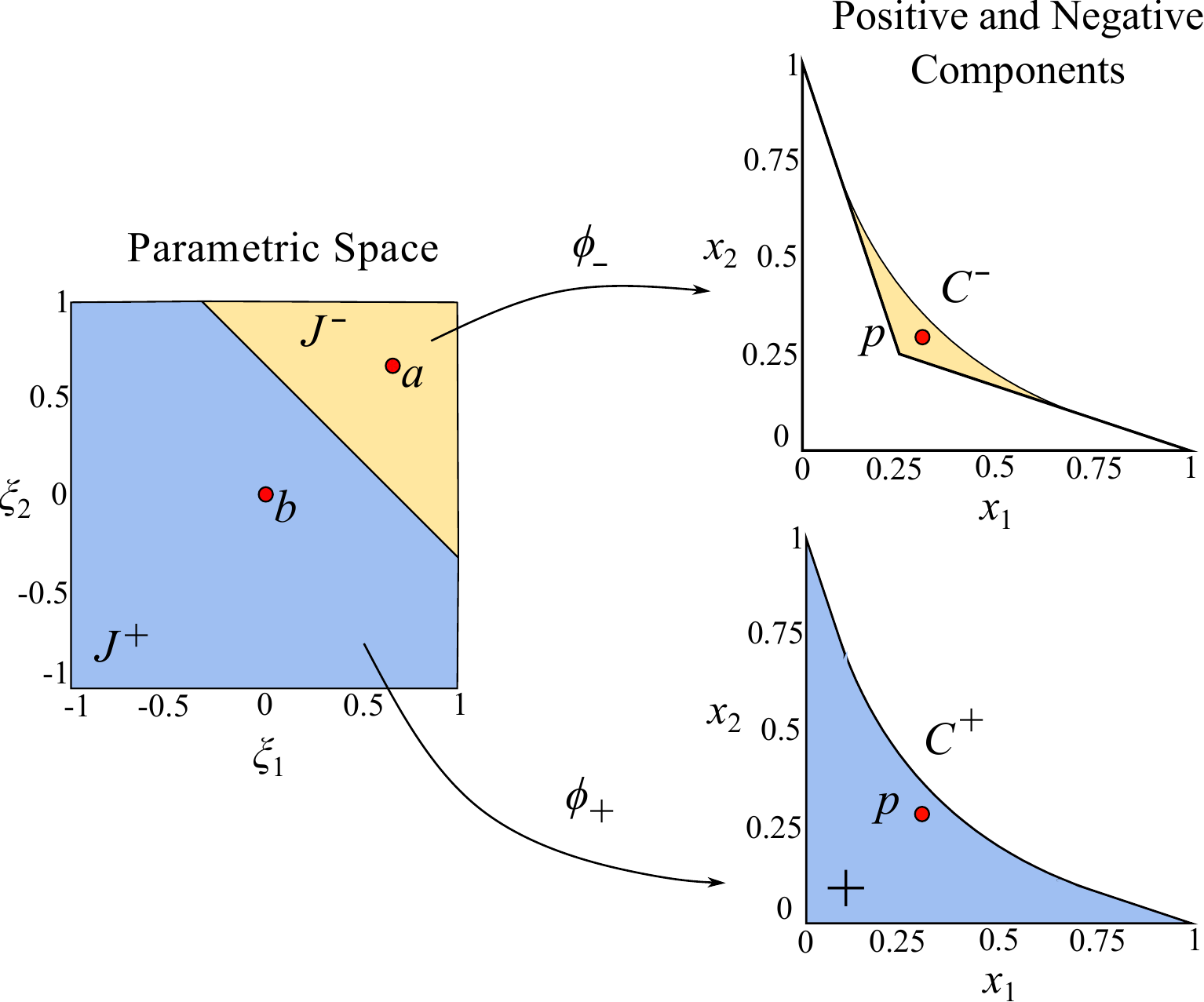}
		\caption{}
		\label{Fig_ConcaveQualMap_d}
	\end{subfigure}	
	\caption{(a) Parametric space of the concave element. Parametric space can be divided into positive and negative Jacobian regions. (b) Physical space of the concave Q9 element.  (c) Positive and negative components}
	\label{Fig_ConcaveQualMap}	
\end{figure}

Consider the standard isoparametric mapping from $(\xi_1,\xi_2)$ space  in Fig.~\ref{Fig_ConcaveQualMap_a} to a concave element in the physical space	$(x_1, x_2)$ in  Fig.~\ref{Fig_ConcaveQualMap_b}.  Observe that the element folds onto itself. Further, the parametric space can be divided into positive ($J^+$) and negative ($J^-$) Jacobian regions and the parametric mapping $\phi$ is \emph{not fully invertible}. 

The main idea in i-TFEM is that the positive and negative parametric regions are treated separately, thus relaxing the constraint of full invertibility to piecewise invertibility. In particular, the physical space corresponding to the positive (negative) parametric region $J^+$  ($J^-$) is termed as positive (negative) component and is denoted by $C^+$ ($C^-$). 
Observe the piecewise mapping
$$
\phi_\pm \colon J^\pm \rightarrow C^\pm
$$
is invertible i.e. bijective (see Fig.~\ref{Fig_ConcaveQualMap_d}).

Let $\boldsymbol{N}_j(\boldsymbol{\xi})$ be the standard bilinear Lagrange shape functions defined over the parametric space of element $E_j$.   
Let $\boldsymbol{N}_j^\pm$ be the restriction of $\boldsymbol{N}_j$ to $J^\pm$, i.e.,
\begin{equation}
	\boldsymbol{N}_j^\pm(\boldsymbol{x}) \coloneqq  \boldsymbol{N}_j(\phi^{-1}_\pm(\boldsymbol{x}))
\end{equation}
The corresponding field is then given by
\begin{equation}
	\boldsymbol{u}_j^\pm(\boldsymbol{x}) = \boldsymbol{N}_j^\pm(\boldsymbol{x})  \boldsymbol{{\hat u}}_j
\end{equation}

\begin{figure}[H]
	\centering\includegraphics[width=1\linewidth]{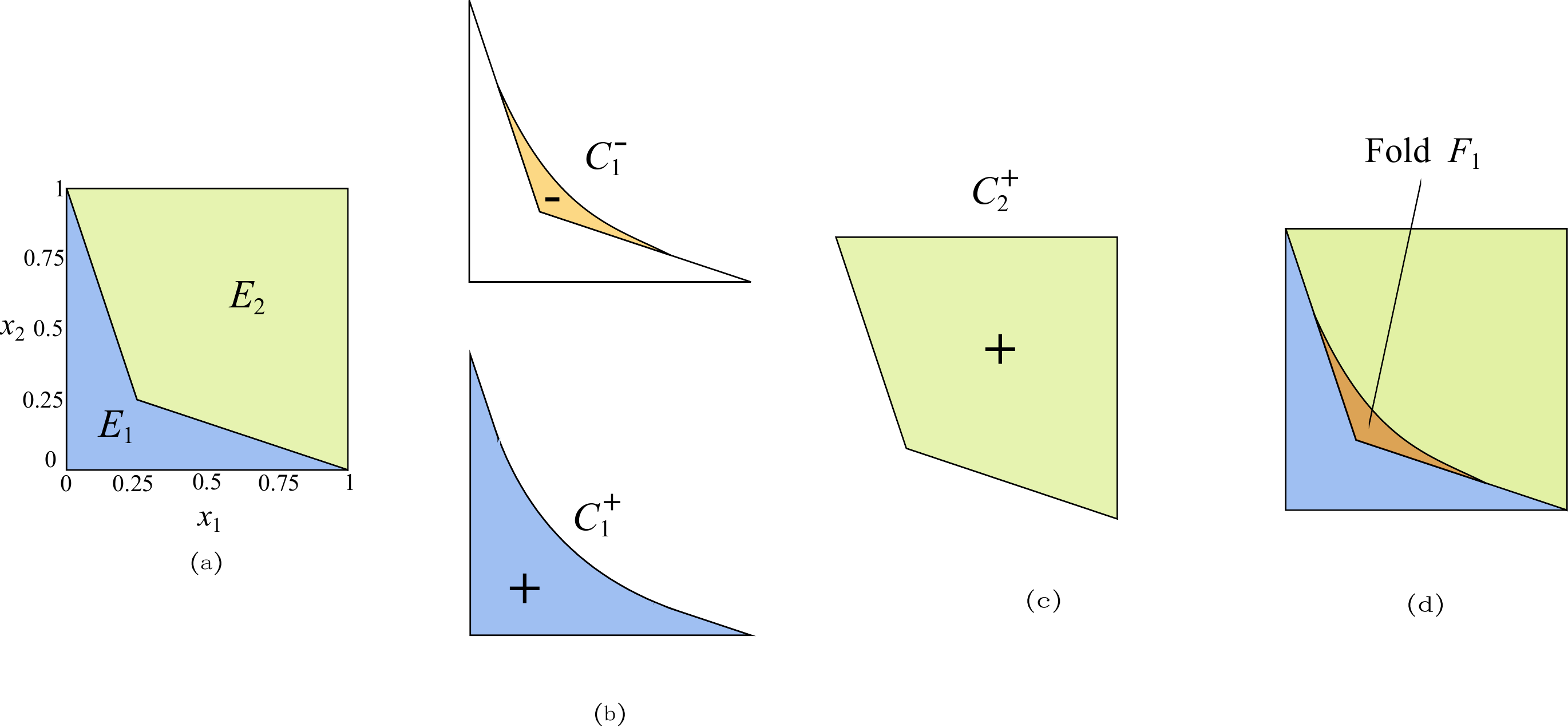}
	\caption{(a) 2-D domain discretized into two bilinear quads. (b) Positive and negative $|\boldsymbol{J}|$ regions of the concave element. (c) Convex element of the mesh. (d) Final physical space is self-overlapping. }
	\label{Fig_PositiveNegativeJRegions}
\end{figure}

\begin{figure}[H]
	\centering\includegraphics[width=0.65\linewidth]{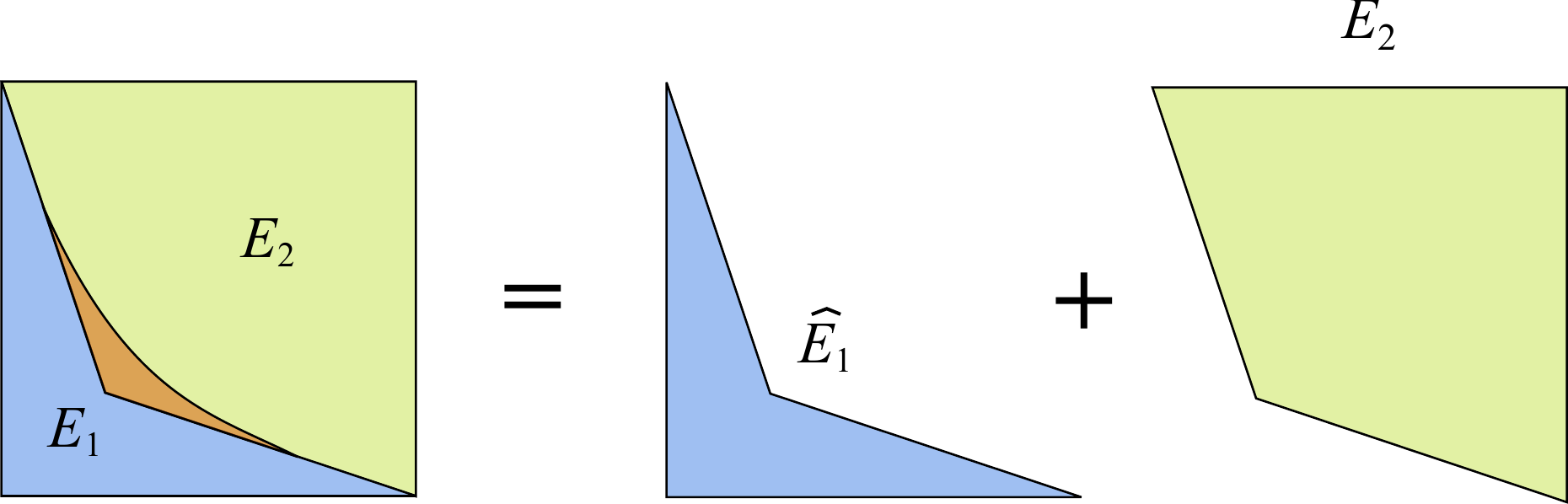}
	\caption{Parts contributing to the field definition.}
	\label{Fig_C1divisions}
\end{figure}

Now consider the two-element tangled patch in Fig.~\ref{Fig_PositiveNegativeJRegions}a. The positive and negative components ($C_1^+$ and $C_1^-$) of the concave element $E_1$ are shown in Fig.~\ref{Fig_PositiveNegativeJRegions}b.  On the other hand, the convex element $E_2$ has only one positive component (see Fig.~\ref{Fig_PositiveNegativeJRegions}c):	$E_2 = C_2^+$ while  $ C_2^- = \emptyset$. Further, the fold $F_1$ illustrated in Fig.~\ref{Fig_PositiveNegativeJRegions}d overlaps with $E_2$ as well.
Thus, for any point $\boldsymbol{x} \in F_1$, all the three components overlap, and three fields  $\boldsymbol{u}^{+}_1\left( \boldsymbol{x}\right)$, $\boldsymbol{u}^{-}_1\left( \boldsymbol{x}\right)$ and $\boldsymbol{u}_2^+(\boldsymbol{x})$ can be defined. 
\begin{comment}
	Further, since $\boldsymbol{x}$ also belongs to the element $E_2$, one can define the corresponding shape functions $\boldsymbol{N}_2^+(\boldsymbol{x})$. 
These shape functions are used to compute the approximate fields over the fold as:
\begin{equation}
	\boldsymbol{u}_1^+(\boldsymbol{x}) \coloneqq \boldsymbol{N}_1^+(\boldsymbol{x})  \boldsymbol{{\hat u}}_1, \quad \boldsymbol{u}_1^-(\boldsymbol{x}) \coloneqq \boldsymbol{N}_1^-(\boldsymbol{x})  \boldsymbol{{\hat u}}_1, \quad \text{and}  \quad \boldsymbol{u}_2^+(\boldsymbol{x}) \coloneqq \boldsymbol{N}_2^+(\boldsymbol{x})  \boldsymbol{{\hat u}}_2
\end{equation}
\end{comment}
 Thus, the field is clearly ambiguous within the fold. 
 Removing the ambiguity in the field definition is the first step in i-TFEM \cite{prabhune2023computationally, prabhune2022tangled, danczyk2013finite}. In particular, in i-TFEM, we \emph {define} the field at a point $\boldsymbol{x}$ within a fold as:
 \begin{equation}
 	\boldsymbol{u}(\boldsymbol{x}) \coloneqq \boldsymbol{u}_2^+(\boldsymbol{x}), \quad \forall \boldsymbol{x} \in F_1 
 	\label{Eq_fieldAtqReduced}
 \end{equation}
The underlying reasons are discussed in \cite{prabhune2023computationally}, but briefly, this is necessary for field continuity and to capture constant strain fields. In other words, for iso-parametric elements, the tangled region can be considered as being part of just the convex element $E_2$. This naturally leads to a division of the mesh into two parts: $E_2$ and $\widehat{E}_1$ as illustrated in Fig.~\ref{Fig_C1divisions}. Note that $\widehat{E}_1$ does not include the folded region whereas $E_1$ does (this is elaborated further below). 

Thus, the field $\boldsymbol{u}$ over the two-element patch is defined as:
\begin{equation}
	\boldsymbol{u}(\boldsymbol{x}) |_{\widehat{E}_1} = \boldsymbol{N}^+_1(\boldsymbol{x}) \boldsymbol{\hat{u}}_1 \quad \text{and} \quad \boldsymbol{u}(\boldsymbol{x}) |_{{E}_2} = \boldsymbol{N}^+_2 (\boldsymbol{x}) \boldsymbol{\hat{u}}_2
	\label{Eq_finalFieldDef}
\end{equation}

However, an additional constraint is needed to ensure the continuity of the field across their common boundary. In particular, by approaching the re-entrant corner, from $\widehat{E}_j$ and $E_2$, one can show that field continuity across the entire boundary is satisfied if and only if (see \cite{prabhune2023computationally}):
\begin{equation}
	\boldsymbol{u}_1^+(\boldsymbol{x}) -  \boldsymbol{u}_1^-(\boldsymbol{x}) = 0, \quad \boldsymbol{x} \in F_1
\end{equation}
This not only makes the field continuous, it also  and forces the contribution of the concave element to be zero in the folded region. 

Introducing the notation $\llbracket     \cdot \rrbracket  = (\cdot)^+ -  (\cdot)^-  $,  the above constraint can be written as:
\begin{equation}
	\llbracket      \boldsymbol{u} _ 1 \rrbracket = 0, \quad \text{in } F_1 
	\label{Eq_equalityCondition}
\end{equation}
In summary, for any pair of overlapping elements $E_1$ and $E_2$ (1) we decompose them into fully invertible regions $\widehat{E}_1$ and $E_2$,
 %(2) for stiffness matrix computations, only the parametric region corresponding to $\widehat{E}_1$ is considered, and (3)
and (2) the constraint (Eq. \ref {Eq_equalityCondition}) is enforced. 
 
\subsection{Variational formulation}
We now consider the implications of these concepts in non-linear elasticity. Our objective is to generalize the residual in Eq. \ref {Eq_R_FEM} and the iteration in Eq.  \ref {Eq_FEMFiniteStrain}, to account for tangling.  Towards this end, we modify the potential energy functional as follows: 
\begin{equation}
	\begin{split}
	\widetilde \Pi =  	\sum_{j \in I_\text{convex}} \; \int\limits_{E_{j}} 	\left( \Psi	\left(\boldsymbol{u}_j \right)  -			\boldsymbol{u}	_j \cdot	\boldsymbol{b}\right) 	dV + 				\sum_{j \in I_\text{concave}} \;  \int\limits_{\widehat{E}_{j}}	\left( \Psi	\left(\boldsymbol{u}_j \right)  -			\boldsymbol{u}_j	\cdot	\boldsymbol{b}\right) 	dV -  \sum_{j \in I}	 \int\limits_{\partial E_{j}^T}	{\boldsymbol{u}_j \cdot	\boldsymbol{T}dS} \\ 	    + \sum_{j \in I_\text{concave}} \; \int\limits_{F_j}	{	\boldsymbol{\lambda}_j \cdot 	 \llbracket     \boldsymbol{u} _ j\rrbracket 	 }dV
	\end{split}
\label{Eq_modifiedPi_2elem}
\end{equation}
where the concave and convex elements are indexed as $I_\text{concave}$ and $I_\text{convex}$ respectively, and  the constraints in Eq. \ref {Eq_equalityCondition} are included via Lagrange multipliers $\boldsymbol{\lambda} $.

We now set the variation of the potential energy with respect to $\boldsymbol{u}$ and $\boldsymbol{\lambda}$ to zero:
\begin{equation}
	\delta_{\boldsymbol{u},\boldsymbol{\lambda}} \widetilde \Pi = 0.
\end{equation}
This leads to the following weak form where $\forall \delta \boldsymbol{u} \in H^1_0$ and $\forall \delta \boldsymbol{\lambda} \in L^2$, i.e., we seek  $\boldsymbol{u} \in H^1$ and $\boldsymbol{\lambda} \in L^2$ such that: 
\begin{equation}
	\begin{split}
			\left[  \sum_{j \in I_\text{convex}} \int\limits_{E_j}  	\left(  \boldsymbol{P} : \frac{\partial	\delta	\boldsymbol{u}_j}	{\partial\boldsymbol{X}} -	\delta\boldsymbol{u}_j	\cdot	\boldsymbol{b}\right) 	dV 																							+																				 \sum_{j \in I_\text{concave}} \int\limits_{\widehat{E}_j}  	\left(  \boldsymbol{P} : \frac{\partial	\delta	\boldsymbol{u}_j}	{\partial\boldsymbol{X}} -	\delta\boldsymbol{u}_j	\cdot	\boldsymbol{b}\right) 	dV  																					- 																					\sum_{j \in I} \int\limits_{\partial	E_j ^T }	{\delta\boldsymbol{u}_j	\cdot	\boldsymbol{T}dS} 	 \right. \\
		\left. +																						  \sum_{j \in I}  \int\limits_{F_j}{ \boldsymbol{\lambda}_j  \cdot \llbracket    \delta \boldsymbol{u} _ j\rrbracket    }dV	\right] 		 +		 					
		 \int\limits_{F_j}{	\delta\boldsymbol{\lambda}_j \cdot	 \llbracket     \boldsymbol{u} _ j\rrbracket   }dV								  = 0
	\end{split}
\label{Eq_weakForm}
\end{equation}
where $\boldsymbol{P}$ is the first Piola-Kirchhoff stress tensor. Next, we approximate the fields using the standard (Bubnov-) Galerkin formulation:
\begin{equation}
\boldsymbol{u}_j  = \boldsymbol{N}_j \boldsymbol{{\hat u}}_j,  \quad  \boldsymbol{\lambda} _j = \boldsymbol{N}^\lambda _j \boldsymbol{\hat{\lambda}}_j
\label{Eq_fieldApprox}
\end{equation}
This leads to:
\begin{equation}
\delta\boldsymbol{\hat{u}}^\top \widetilde{\boldsymbol{R}} + \delta\boldsymbol{\hat{\lambda}}^\top \boldsymbol{{C}}^\top \boldsymbol{\hat{u}} = 0
\label{Eq_FiniteDimWeakForm}
\end{equation}
Here, $\delta\boldsymbol{\hat{u}}^\top \widetilde{\boldsymbol{R}}$ represents the terms in the square bracket of Eq.~\ref{Eq_weakForm}. Observe from Eq. \ref{Eq_weakForm}, that the residual $\widetilde{\boldsymbol{R}}  $ can be expressed as: 
\begin{equation}
	\widetilde{\boldsymbol{R}} \left( \boldsymbol{\hat{u}}, \boldsymbol{\hat{\lambda}} \right)  = \boldsymbol{R}^{u}\left( \boldsymbol{\hat{u}} \right)  + \boldsymbol{C} \boldsymbol{\hat{\lambda}} = \boldsymbol{0}
	\label{Eq_R_iTFEM}
\end{equation}
where $\boldsymbol{R}^{u}$ involves only the primary field  $\boldsymbol{u}$ and requires integrating over both the convex and concave elements:
 \begin{equation}
	\boldsymbol{R}^u = \boldsymbol{R}^u_{convex} +   \widehat{\boldsymbol{R}}^u_{concave}
	\label{Eq_Rt_isoparametric}	
\end{equation}
As one can easily deduce, the computation of $\boldsymbol{R}^u_{convex}$ is as in standard FEM. However, the integration over the concave elements must be carried out over the subset of parametric space; see Fig.~\ref{Fig_parametricSpace}. This is discussed in detail in the next subsection. 

Next, to solve Eq. \ref{Eq_R_iTFEM} through iterations, we consider the first order Taylor series:
\begin{equation}
	\boldsymbol{R}^u\left(\boldsymbol{\hat{u}}^n\right) 		+ 			\frac{\partial \boldsymbol{R}^u}{\partial \boldsymbol{\hat{u}}} \bigg{|}_n \Delta \boldsymbol{\hat{u}}^{n+1} 			+
	\boldsymbol{C} \boldsymbol{\hat{\lambda}}^n +
	 \boldsymbol{C} \Delta \boldsymbol{\hat{\lambda}}^{n+1}	 	= 		\boldsymbol{0}
	\label{Eq_RNewtonRaphson}
\end{equation}
i.e.,
\begin{equation}
	\boldsymbol{K}^{t}  \Delta \boldsymbol{\hat{u}} +  \boldsymbol{C} \Delta\boldsymbol{\hat{\lambda}}= - \left( \boldsymbol{R}^u +   \boldsymbol{C} \boldsymbol{\hat{\lambda}}\right) 
	\label{Eq_NREq1} 
\end{equation}
where
\begin{equation}
	\boldsymbol{K}^{t}   = \boldsymbol{K}_{convex}^{t} + \widehat{\boldsymbol{K}}_{concave}^{t}
	\label{Eq_KNL} 
\end{equation}
Here,  $\boldsymbol{K}^{t}_{convex}$ and $\widehat{\boldsymbol{K}}_{concave}^{t}$ are tangent matrices corresponding to convex and concave elements respectively. 
Further, from Eq. \ref {Eq_FiniteDimWeakForm}, we have
\begin{equation}
	\boldsymbol{C}^\top \Delta \boldsymbol{\hat{u}}  = \boldsymbol{0}
	\label{Eq_NREq2} 
\end{equation}
From Eq.~\ref{Eq_NREq1} and Eq.~\ref{Eq_NREq2} we have the final set of linear equations one must solve iteratively:
\begin{equation}
	\begin{bmatrix}
		\boldsymbol{K}^{t} & \boldsymbol{C}\\
		\boldsymbol{C}^\top& \boldsymbol{0}
	\end{bmatrix} 	
	\begin{Bmatrix}
		\Delta \boldsymbol{\hat{u}}^{n+1}\\
		\Delta\boldsymbol{\hat{\lambda}}^{n+1}
	\end{Bmatrix} 	= 
	\begin{Bmatrix}
		-\left( \boldsymbol{R}^u + \boldsymbol{C} \boldsymbol{\hat{\lambda}} \right) \\
		\boldsymbol{0}
	\end{Bmatrix}.
	\label{Eq_GeneralizedimplicitTFEMFiniteStrain}
\end{equation}
If the mesh does not contain any tangled elements, then $\boldsymbol{K}^{t}   = \boldsymbol{K}_{convex}^{t}$, and $ \boldsymbol{C} $ does not exist, i.e., i-TFEM reduces to standard FEM.

\subsection{Implementation Details}
We now discuss the implementation considering a two-element mesh. 
\subsubsection{Computation of residual vector and stiffness matrix}
As mentioned earlier, to compute $\boldsymbol{R}^u_{convex}$ in Eq. \ref {Eq_Rt_isoparametric}, standard FEM procedures with Gauss quadrature can be used.  However, to compute  $\boldsymbol{R}^u_{concave}$, only the fully invertible subset $\widehat{E}_1$ is to be considered.  

\begin{figure}[H]
	\centering\includegraphics[width=0.7\linewidth]{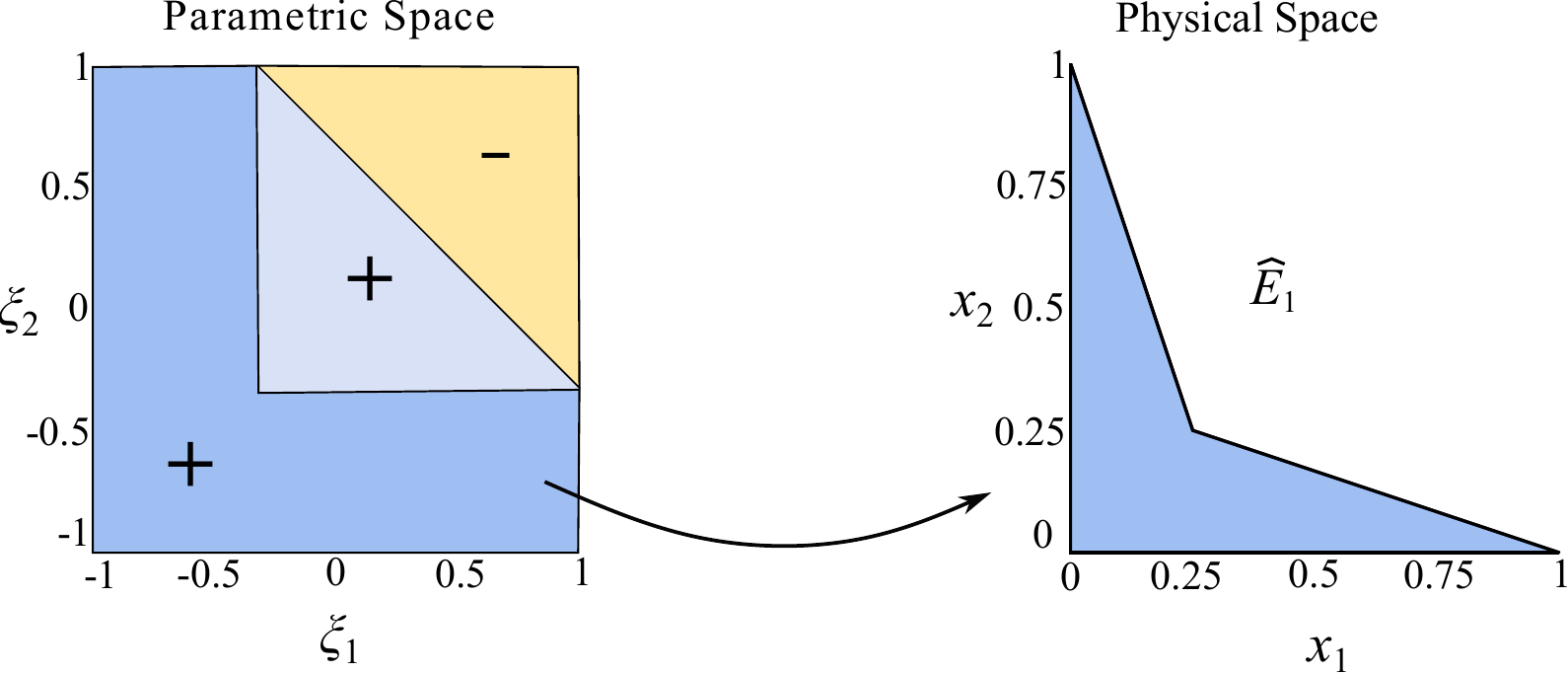}
	\caption{(a) Parametric space. (b) $\widehat{E}_1$}
	\label{Fig_parametricSpace}
\end{figure}

Observe that $\widehat{E}_1$ is not the same as $E_1$. Though both have the same physical boundary, they represent different regions of parametric space. Specifically, $E_1$ represents the entire parametric space, while $\widehat{E}_1$ represents only the L-shaped subset of the positive $ |\boldsymbol{J}|$ region as illustrated in Fig.~\ref{Fig_parametricSpace}.  In other words, $ \widehat{E}_1$  corresponds to  a \emph{fully invertible} subset of the parametric space of the concave element.

\begin{figure}[H]
	\centering\includegraphics[width=0.2\linewidth]{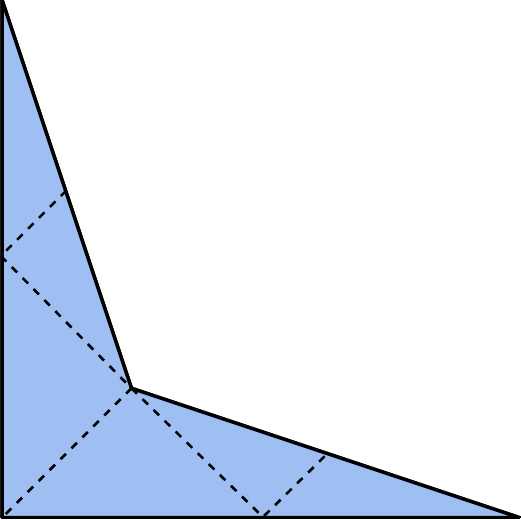}
	\caption{Triangulation of  $\widehat{E}_1$}
	\label{Fig_triangulatedRegion_concaveElemDotted}
\end{figure}

Therefore, standard Gauss quadrature cannot be employed; instead, $\widehat{E}_1$ is triangulated as illustrated in Fig.~\ref{Fig_triangulatedRegion_concaveElemDotted}. The triangulation is used merely for the purpose of integration and does not lead to additional  degrees of freedom in i-TFEM.
Similarly, to compute $\boldsymbol{K}_{convex}^{t}$, standard FEM  procedures can be used. However, to compute  $\boldsymbol{K}_{concave}^{t}$, the triangulation in Fig.~\ref{Fig_triangulatedRegion_concaveElemDotted} must be used.
\subsubsection{Constraint enforcement}

Finally, to compute the constraint matrix $\boldsymbol{C}$, note that the interpolation of the Lagrange multiplier $\boldsymbol{\lambda}$ needs to be only square integrable since  its gradient does not appear in the formulation.  For Q4 elements, the primary field $\boldsymbol{u}$ is approximated using standard bilinear functions while $\boldsymbol{N}^\lambda$ are constant functions. Thus, the finite dimensional approximation for the Lagrange multiplier $\boldsymbol{\lambda}$ comes from an FE space that is smaller than that for  $\boldsymbol{u}$. Accordingly, we can write the constraint matrix defined in Equation \ref{Eq_R_iTFEM} as 
\begin{equation}
	\boldsymbol{C} =	\int\limits_{F_{1}} \left( \boldsymbol{N}_1^+  - \boldsymbol{N}_1^- \right)^\top dV.
	\label{Eq_weakEqualityCond2}
\end{equation}
% Observe that the interpolation of the Lagrange multiplier $\boldsymbol{\lambda}$ needs to be only square integrable since we do not consider its gradient in the formulation.

Direct integration over the tangled region $F_1$  to compute $\boldsymbol{C}$ is  computationally expensive and cumbersome \cite{prabhune2022tangled}. Instead, we evaluate the integrand at a sample point $\boldsymbol{x} \in F_1$, say the concave vertex, i.e., evaluate the $\boldsymbol{C}$ as
\begin{equation}
	\boldsymbol{C} =	 \left( \boldsymbol{N}_1^+ (\boldsymbol{p} ) - \boldsymbol{N}_1^- (\boldsymbol{p}) \right)^\top = \llbracket \boldsymbol{N}_1 (\boldsymbol{p}) \rrbracket.
	\label{Eq_weakEqualityCond3}
\end{equation}
where $\boldsymbol{p}$ is the re-entrant vertex. Thus, the constraints can be applied directly as a set of algebraic equations.  Since $\boldsymbol{u}$ is a vector field,  each concave element entails two constraint equations. This is  consistent, for example, with the  algebraic constraints implemented in \cite{prabhune2022tangled, prabhune2022towards, prabhune2023computationally}. 

\section{Numerical Experiments}	
	\label{Section5_Experiments}
In this section, i-TFEM is demonstrated using plane strain nonlinear elasticity problems over various tangled meshes. Numerical experiments are conducted under the following conditions:
	\begin{itemize}			
		\item  The implementation is in MATLAB R2022a, on a standard Windows 10 desktop with Intel(R) Core(TM) i9-9820X CPU running at 3.3 GHz with 16 GB memory. 
		\item  The number of quadrature points for convex quadrilateral elements. 
		\item  The triangulation of a concave element (see Fig.~\ref{Fig_triangulatedRegion_concaveElemDotted}) is performed by employing MATLAB's inbuilt mesher -  \texttt{generateMesh}. The number of quadrature points for triangles is 4. 
		%In 3D, tetrahedralization of concave elements is performed using Tetgen \cite{hang2015tetgen}.  The bounding surfaces are triangulated using  \texttt{generateMesh}, and serve as input to Tetgen. The surface mesh-size is set to a relative size of $h_t$ = 0.05, where $h_t$ is defined as the maximum allowable edge length of a surface triangle. The number of quadrature points for tetrahedrons is chosen to be 4.
		\item The load is applied incrementally in 10 steps. The stopping criteria  for Newton Raphson is $||\Delta \boldsymbol{{\hat u} }|| < 10^{-9} $.
	\end{itemize}	
	Through the experiments, we investigate the following:
	\begin{itemize}
		\item \textbf{Cook's problem, single concave element}: For Cook's membrane problem \cite{cook1974improved},   the error in tip displacement due to the presence of a single concave element is reported as the severity of tangling is increased.
		\item \textbf{Cook's problem, multiple concave elements}: For Cook's membrane problem, with numerous tangled elements: (a)  The displacement at a prescribed location is reported for each load step. (b) Deformed configurations for tangled and regular meshes are also compared. (c) Convergence of the tip displacement as a function of mesh size is studied  and compared against standard FEM. (d) Finally, the convergence rate is evaluated.
		\item \textbf{Punch problem, material non-linearity}: For the punch problem \cite{de2019serendipity}, we include material non-linearity and study the convergence characteristics of i-TFEM.
		\item \textbf{Punch problem, multiple overlap, near-incompressibility}: For the punch problem \cite{de2019serendipity}, a mesh with multiple overlaps is considered for compressible and near-incompressible cases.
		\item \textbf{Thin beam problem}: The performance of i-TFEM in bending dominated response is evaluated with neo-Hookean material model. 
		\item \textbf{Aircraft model}: An example of a tangled mesh is presented to evaluate i-TFEM in practical scenarios.   
	\end{itemize}	
	\subsection{2D Cook's membrane: Single concave element}
\begin{figure}[H]
	\centering\includegraphics[width=0.5\linewidth]{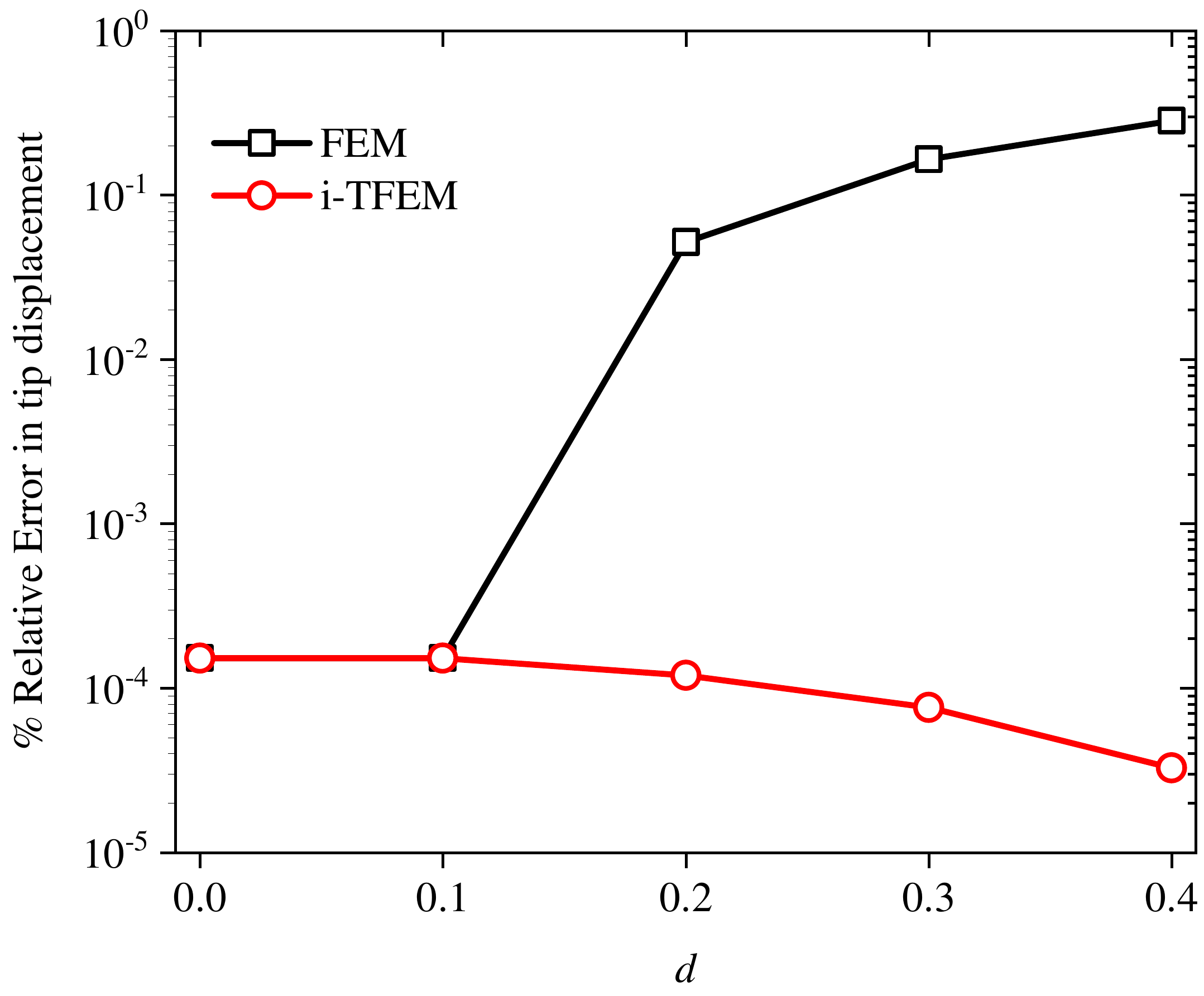}
	\caption{Relative error in tip displacement versus $d$ for FEM and i-TFEM.}
	\label{Fig_cooks1ElemResultsTFEM}
\end{figure}

	To begin with, we solve Cook's membrane problem over the mesh with one concave element illustrated earlier in Fig.~\ref{Fig_cooksMesh1tangle}. Recall that the extent of tangling is controlled by the parameter $d$. For $d >0.1$, a sharp increase in FEM error was observed, as illustrated  in Fig.~\ref{Fig_cooks1ElemResultsTFEM}. On the other hand, using i-TFEM, the error, in fact,  decreases (slightly) for $d >0.1$; see Fig.~\ref{Fig_cooks1ElemResultsTFEM}. 

\subsection{2D Cook's membrane: Multiple concave elements}

\begin{figure}[H]	
		\begin{subfigure}[c]{.5\textwidth}		
		\centering\includegraphics[width=0.6\linewidth]{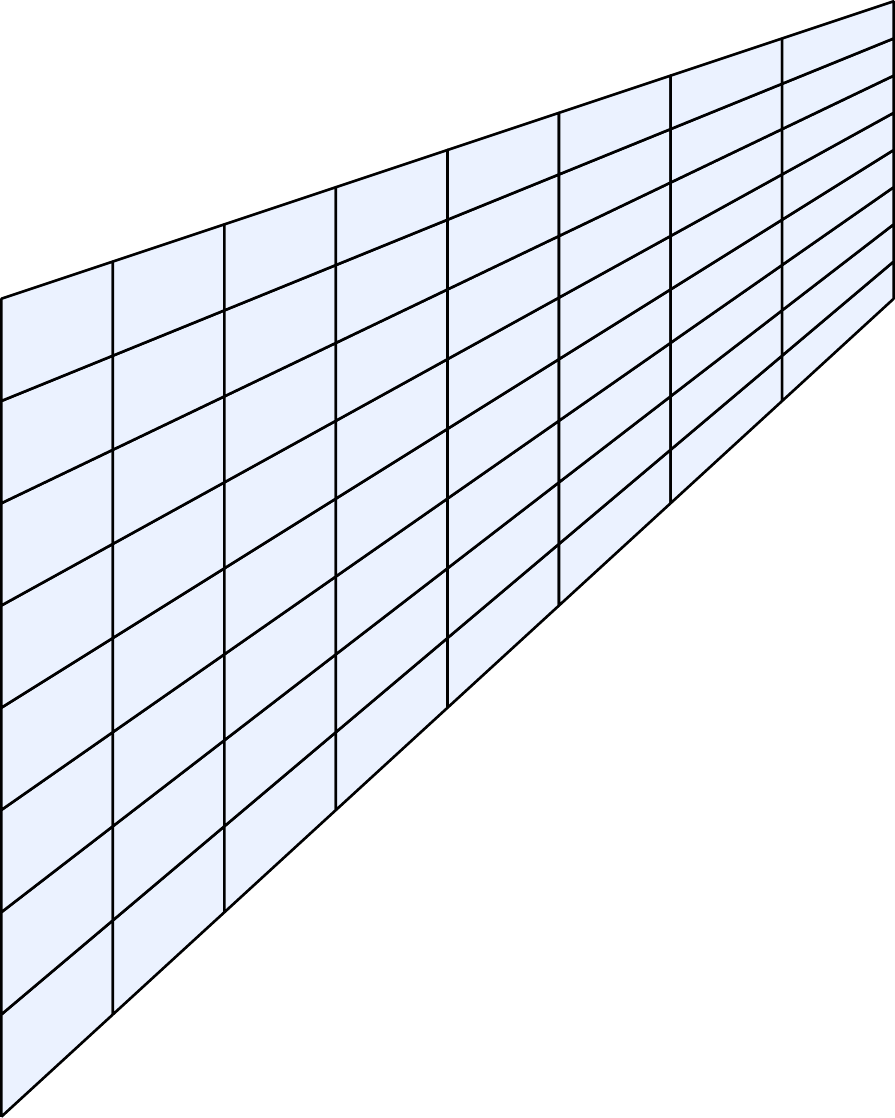}
		\caption{}
		\label{Fig_cooksMeshUntangled2}
	\end{subfigure}	
	\begin{subfigure}[c]{.5\textwidth}			
		\centering\includegraphics[width=0.6\linewidth]{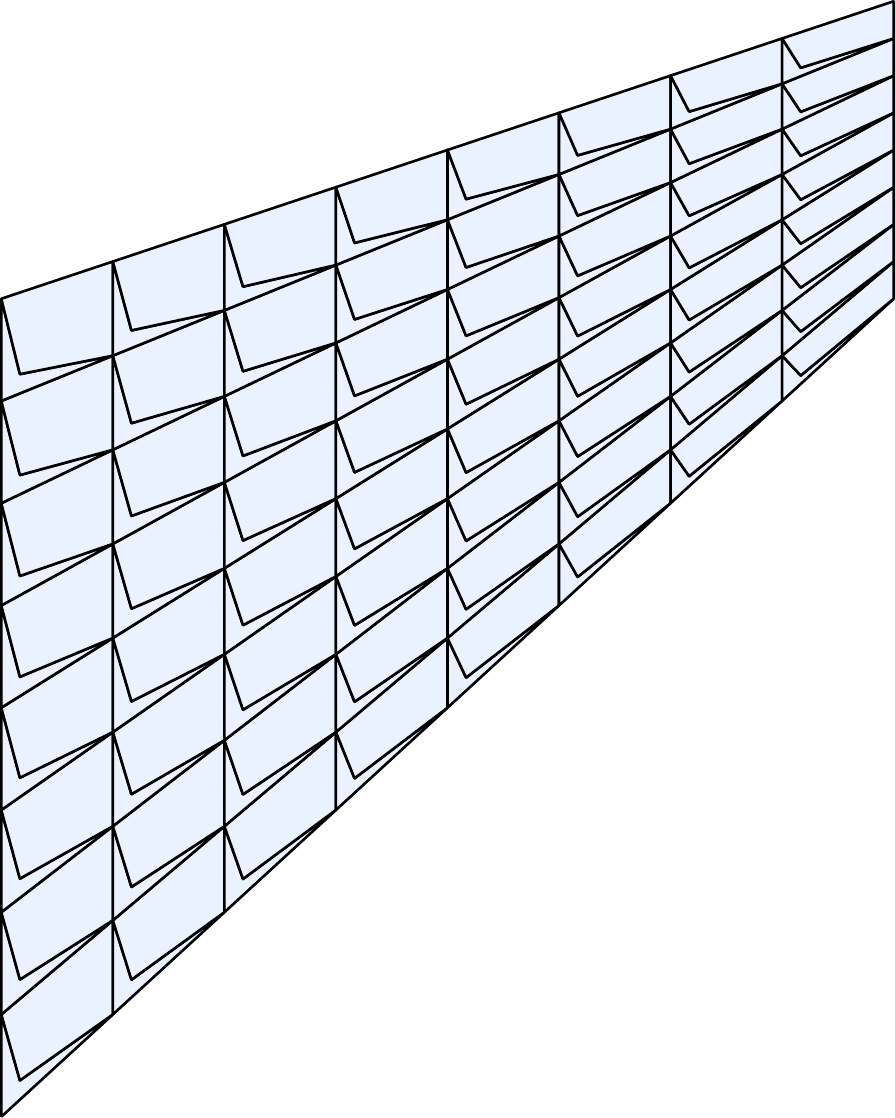}
		\caption{}
		\label{Fig_cooksMeshAllTangle}
	\end{subfigure}

	\caption{Initial configuration  for (a) Regular mesh. (b) tangled mesh with $N = 3 \equiv 8 \times 8$ for the Cook's membrane problem. }
	\label{}	
\end{figure}

Next, we consider a regular mesh illustrated in Fig.~\ref{Fig_cooksMeshUntangled2} and a highly tangled mesh in Fig.~\ref{Fig_cooksMeshAllTangle} where every other element is concave.

\begin{figure}[H]
	\centering\includegraphics[width=0.5\linewidth]{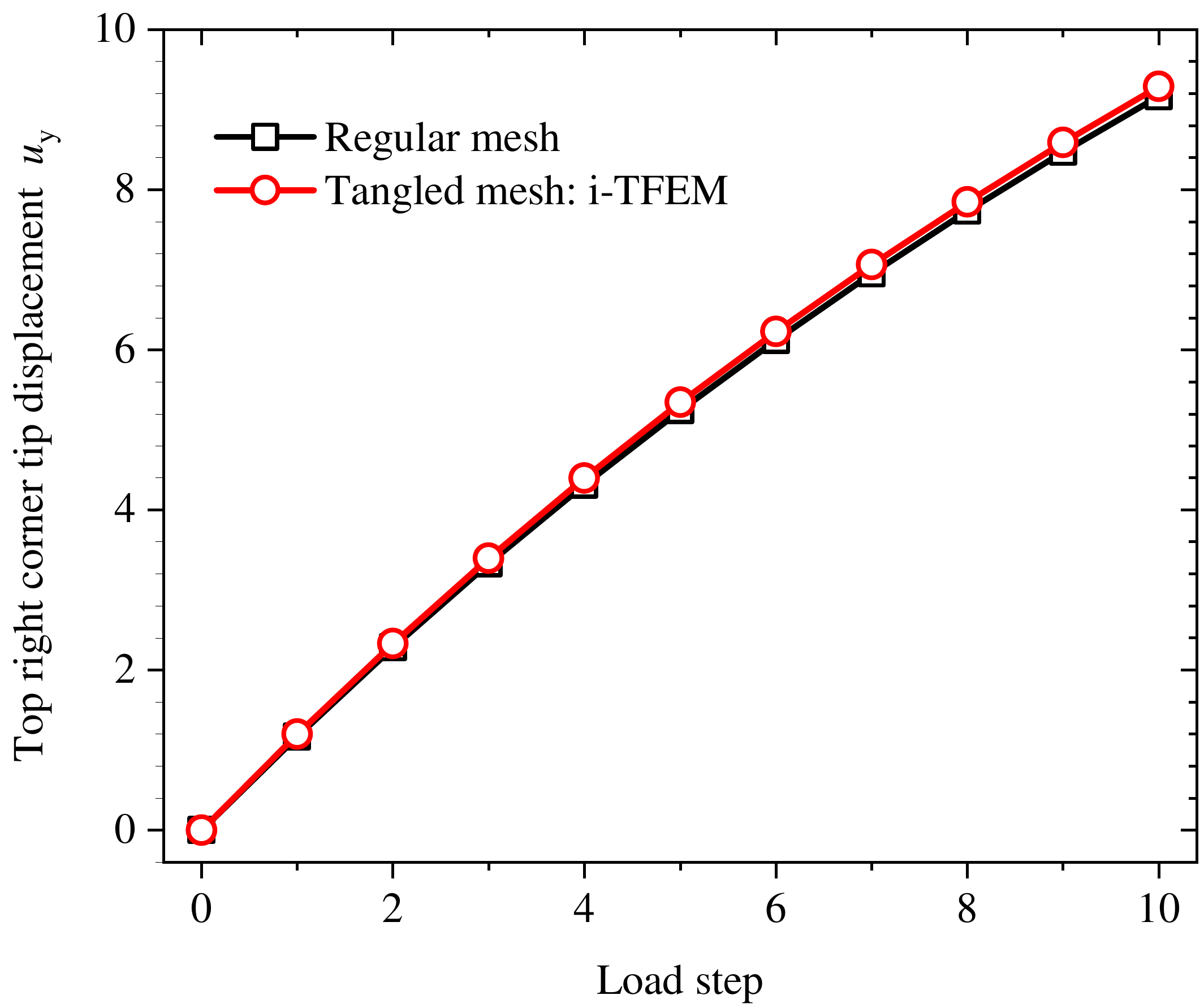}
	\caption{Vertical displacement versus the load step for Cook's membrane problem.}
	\label{Fig_CooksMembraneDispVsLoadStep}
\end{figure}

\begin{figure}[H]
	\begin{subfigure}[c]{.50\textwidth}				
		\centering\includegraphics[width=0.65\linewidth]{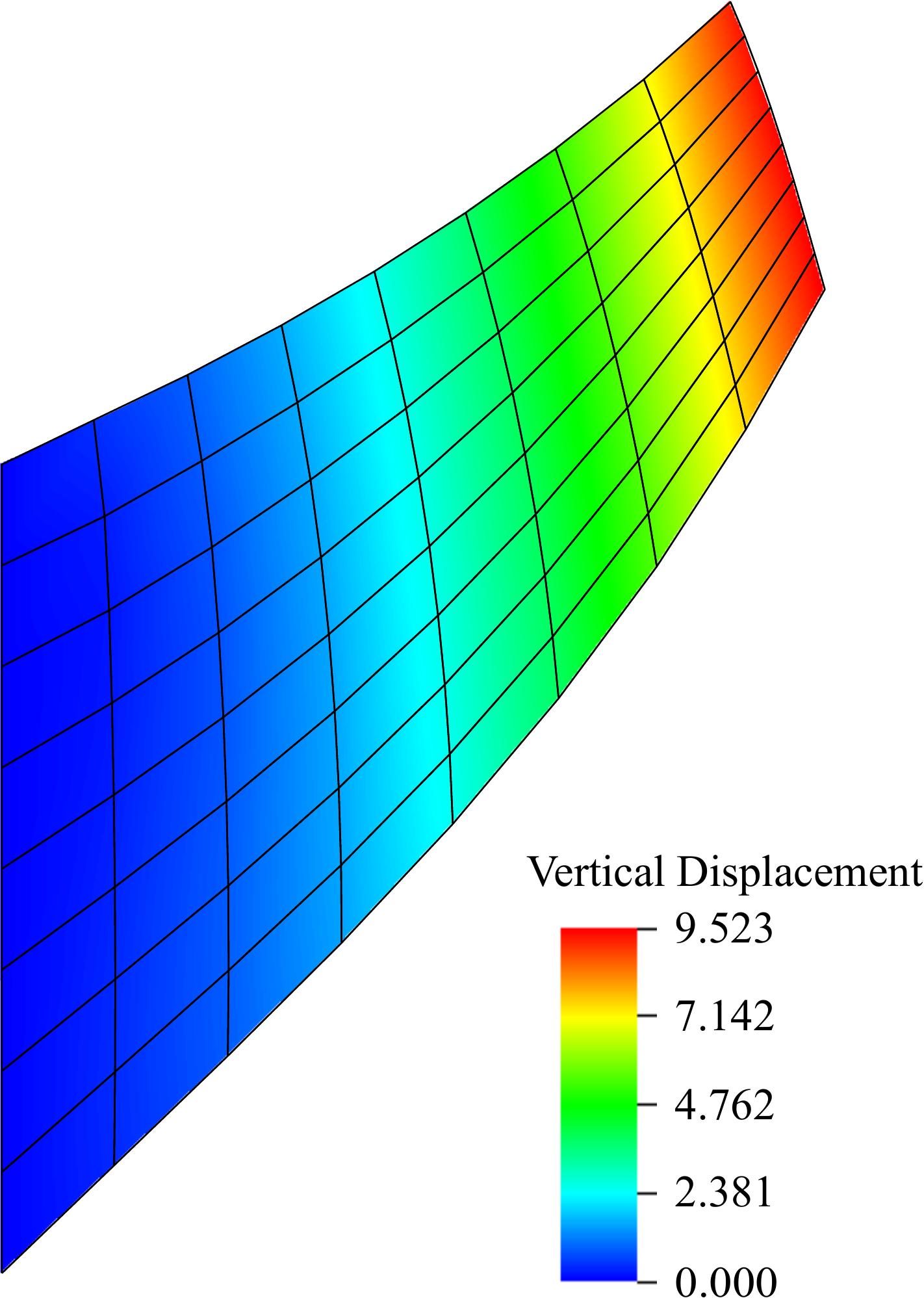}
		\caption{}
		\label{Fig_cooksDisp_untangled}
	\end{subfigure}
	\begin{subfigure}[c]{.50\textwidth}		
		\centering\includegraphics[width=0.65\linewidth]{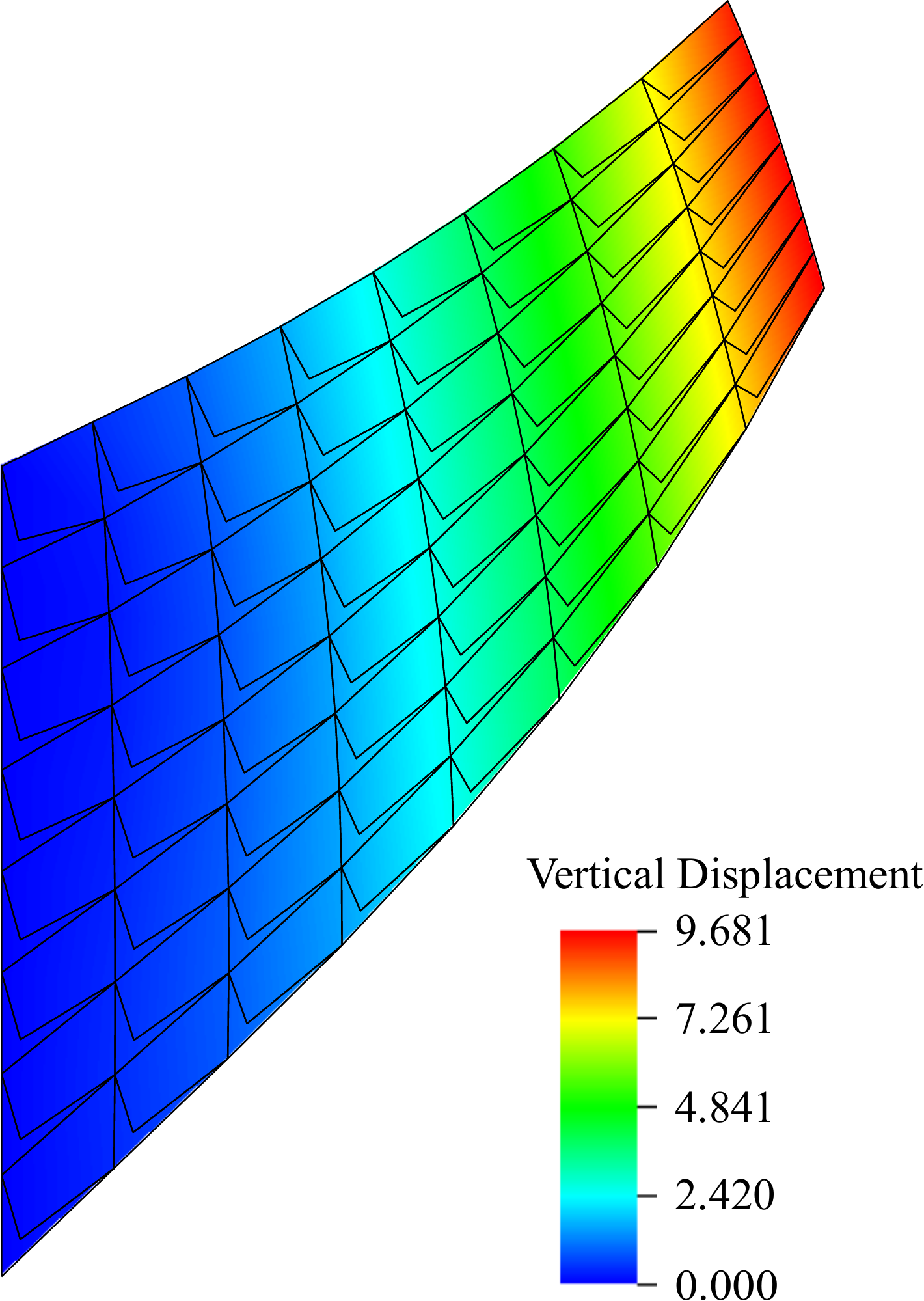}
		\caption{}
		\label{Fig_cooksDisp_TFEM}
	\end{subfigure}	
	\caption{Deformed configuration for (a) Regular mesh and (b) Tangled mesh using i-TFEM for Cook's membrane problem.}
	\label{}	
\end{figure}

The Cook's membrane problem is solved over the regular mesh using standard FEM, and over the tangled mesh using i-TFEM. The vertical displacement at the top right corner point A (see Fig.~\ref{Fig_cooksGeom}) for every load step is reported in Fig.~\ref{Fig_CooksMembraneDispVsLoadStep}. One can observe a  close agreement between the two.
The deformed configuration for regular and tangled meshes after the last load step are reported in Fig.~\ref{Fig_cooksDisp_untangled} and Fig.~\ref{Fig_cooksDisp_TFEM}, respectively. 

\begin{figure}[H]
	\centering\includegraphics[width=0.5\linewidth]{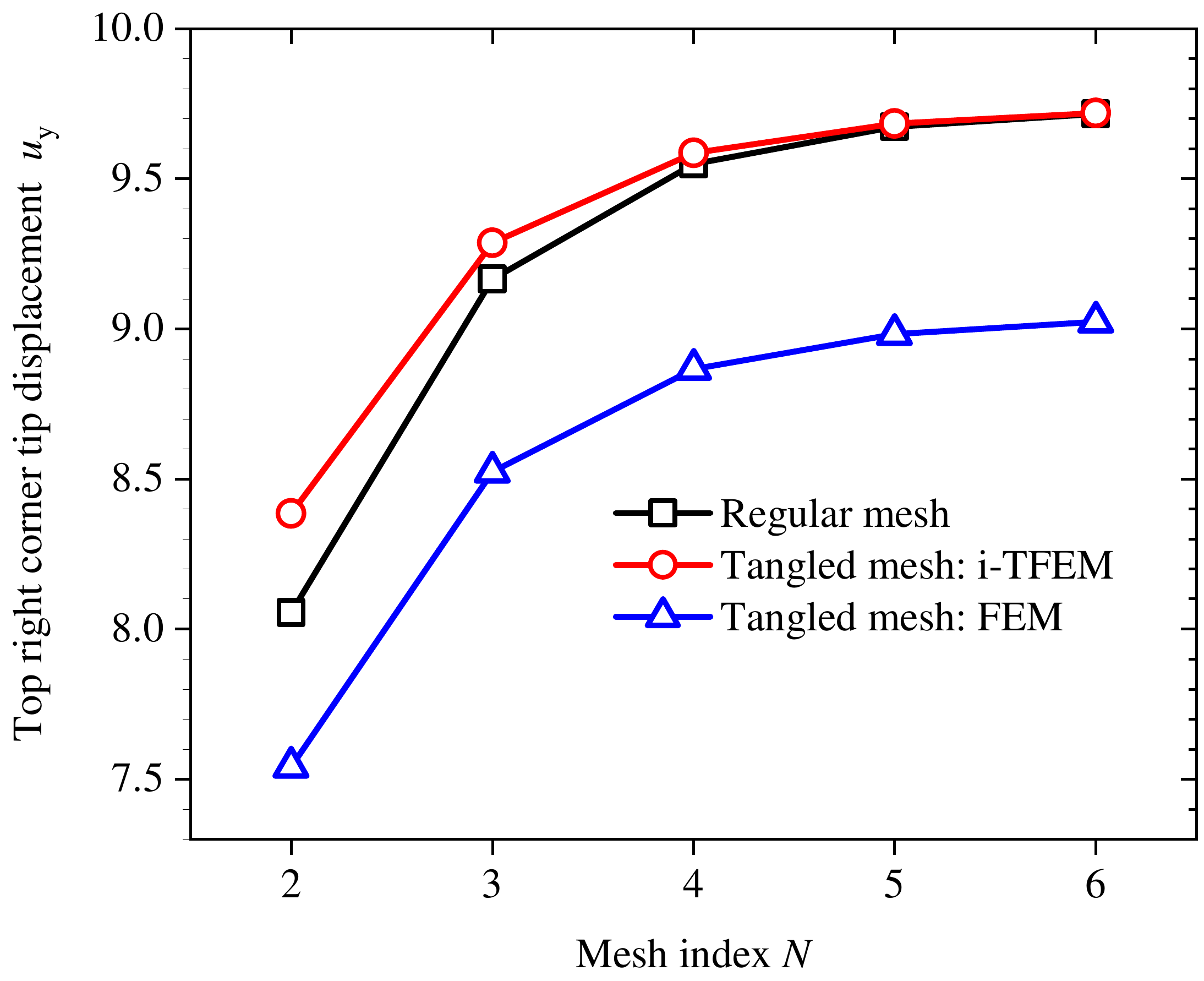}
	\caption{Convergence study for Cook's membrane problem.}
	\label{Fig_CooksMembraneConv}
\end{figure}

\begin{figure}[H]
	\begin{subfigure}[c]{.5\textwidth}				
		\centering\includegraphics[width=0.95\linewidth]{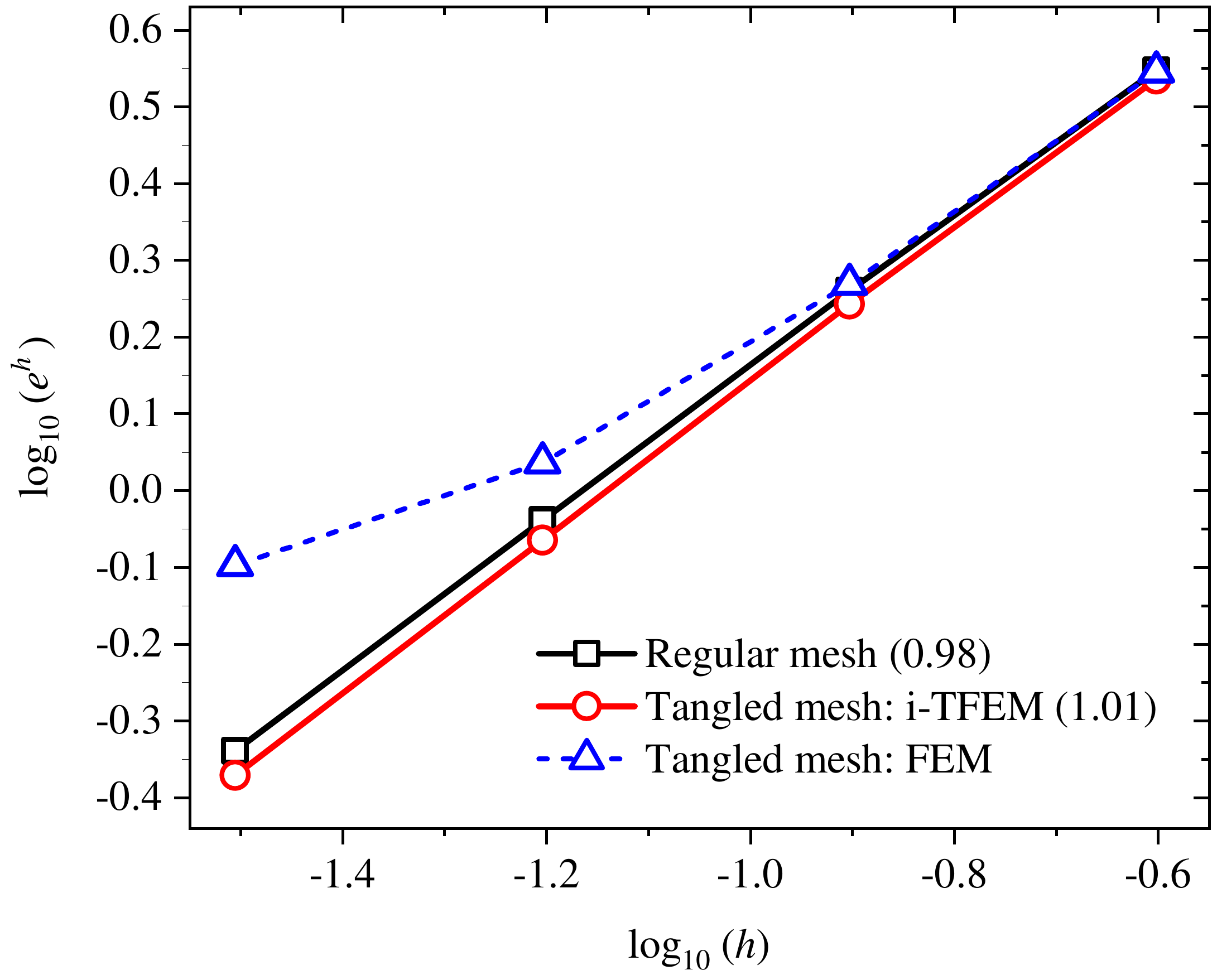}
		\caption{}
		\label{Fig_CooksMembraneConv_H1}
	\end{subfigure}
	\begin{subfigure}[c]{.5\textwidth}				
		\centering\includegraphics[width=0.95\linewidth]{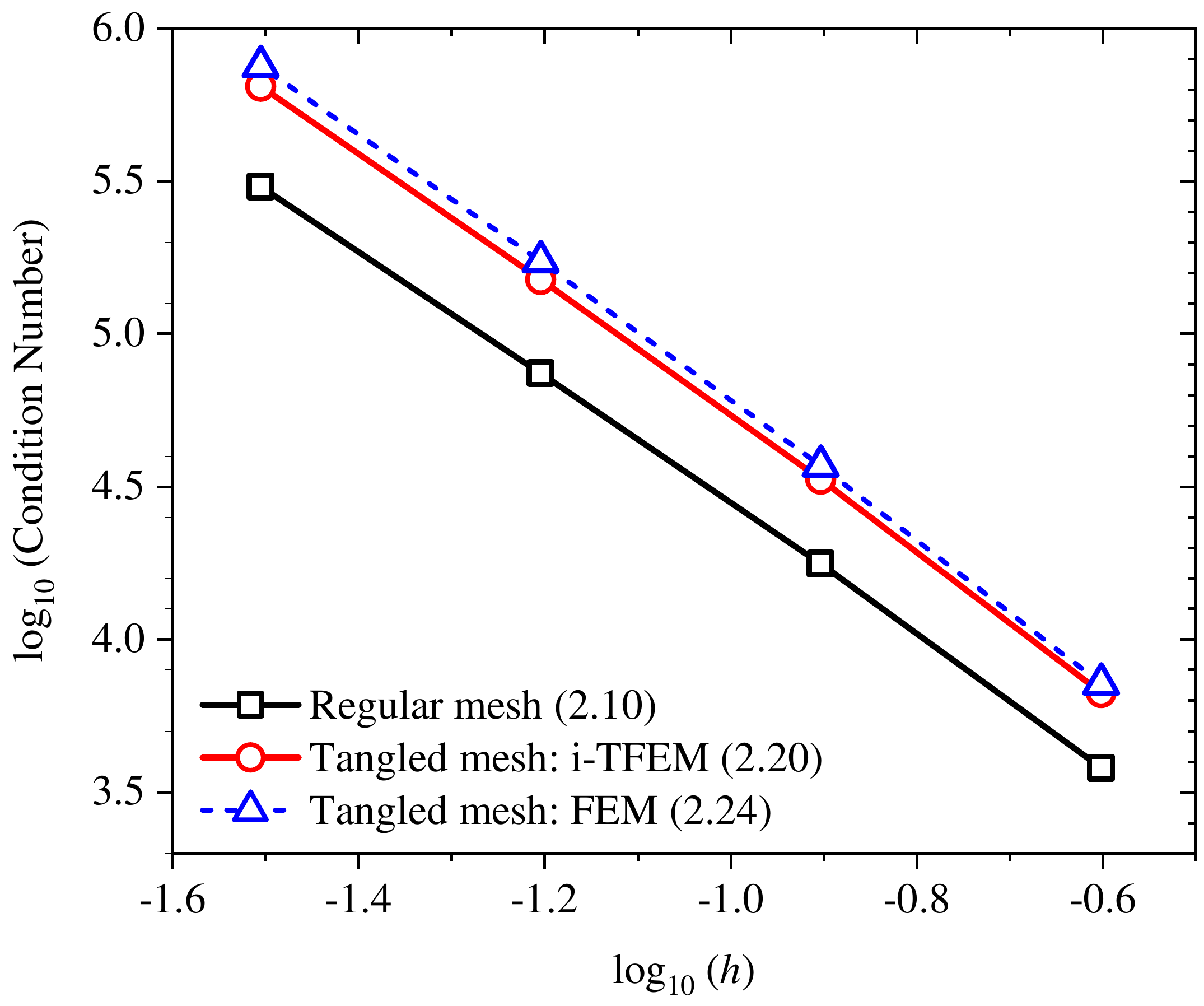}
		\caption{}
		\label{Fig_conditionNum_Cook_h}
	\end{subfigure}	
	
	\caption{(a)${H}^1$ seminorm error versus mesh size and (b) condition number versus mesh size for Cook's problem. The convergence rates are provided in the brackets. }
	\label{}	
\end{figure}

To study convergence, the number of elements is controlled  by a mesh-index $N$, where the number of elements in the regular mesh is $ 2^N \times 2^N$.  Fig.~\ref{Fig_cooksMeshUntangled2} illustrates the regular mesh when $N = 3$, and Fig.~\ref{Fig_cooksMeshAllTangle}, the corresponding tangled mesh. We now compare the solutions from three different methods: standard FEM over regular mesh, standard FEM over tangled mesh, and i-TFEM over tangled mesh. The vertical displacements at point A for all three are plotted as a function of $N$ in Fig.~\ref{Fig_CooksMembraneConv}. Observe that FEM over a regular mesh and i-TFEM over the tangled mesh converge to the same displacement. On the other hand, FEM over a tangled mesh leads to erroneous results.
To study the rate of convergence, we define the $H^1$ seminorm of the displacement error   as  
\begin{equation}
	e^h = {|| \nabla \boldsymbol{u} -  \nabla \boldsymbol{u}^h ||} = \left[ \int\limits_{\Omega} |\nabla \boldsymbol{u} -  \nabla \boldsymbol{u}^h |^2 \; d\Omega \right]  ^{ 0.5} 
\end{equation}	
where $\boldsymbol{u}$ is the reference solution from a fine mesh with $N= 7$, and  $\boldsymbol{u}^h$ is the solution under consideration. Fig.~\ref{Fig_CooksMembraneConv_H1} illustrates the error vs. mesh size ($h$) on a log-log scale over the non-tangled mesh as well as over the tangled mesh using FEM and i-TFEM. One can observe a near-optimal convergence rate for i-TFEM. Next, the effect of mesh size on the condition number of the matrix in Eq.~\ref{Eq_GeneralizedimplicitTFEMFiniteStrain} is studied. 
 Fig.~\ref{Fig_conditionNum_Cook_h} shows that the condition number for tangled meshes increases with mesh size at a rate similar to that of a regular mesh.

\subsection{Punch Problem: Material nonlinearity}
\begin{figure}[H]
	\centering\includegraphics[width=1\linewidth]{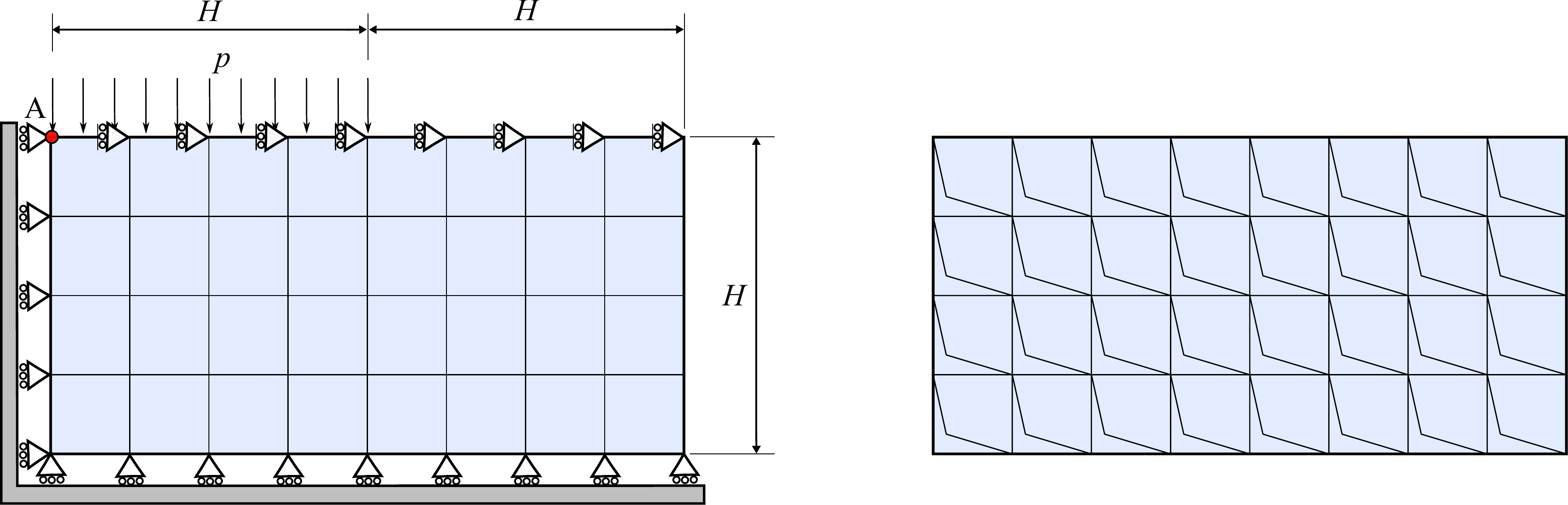}
	\caption{Initial configuration of the punch problem with mesh size $N = 2 \equiv 8 \times 4$ and the tangled  mesh.}
	\label{Fig_punchGeometryNMesh}
\end{figure}

\begin{figure}[H]
 	\centering\includegraphics[width=0.5\linewidth]{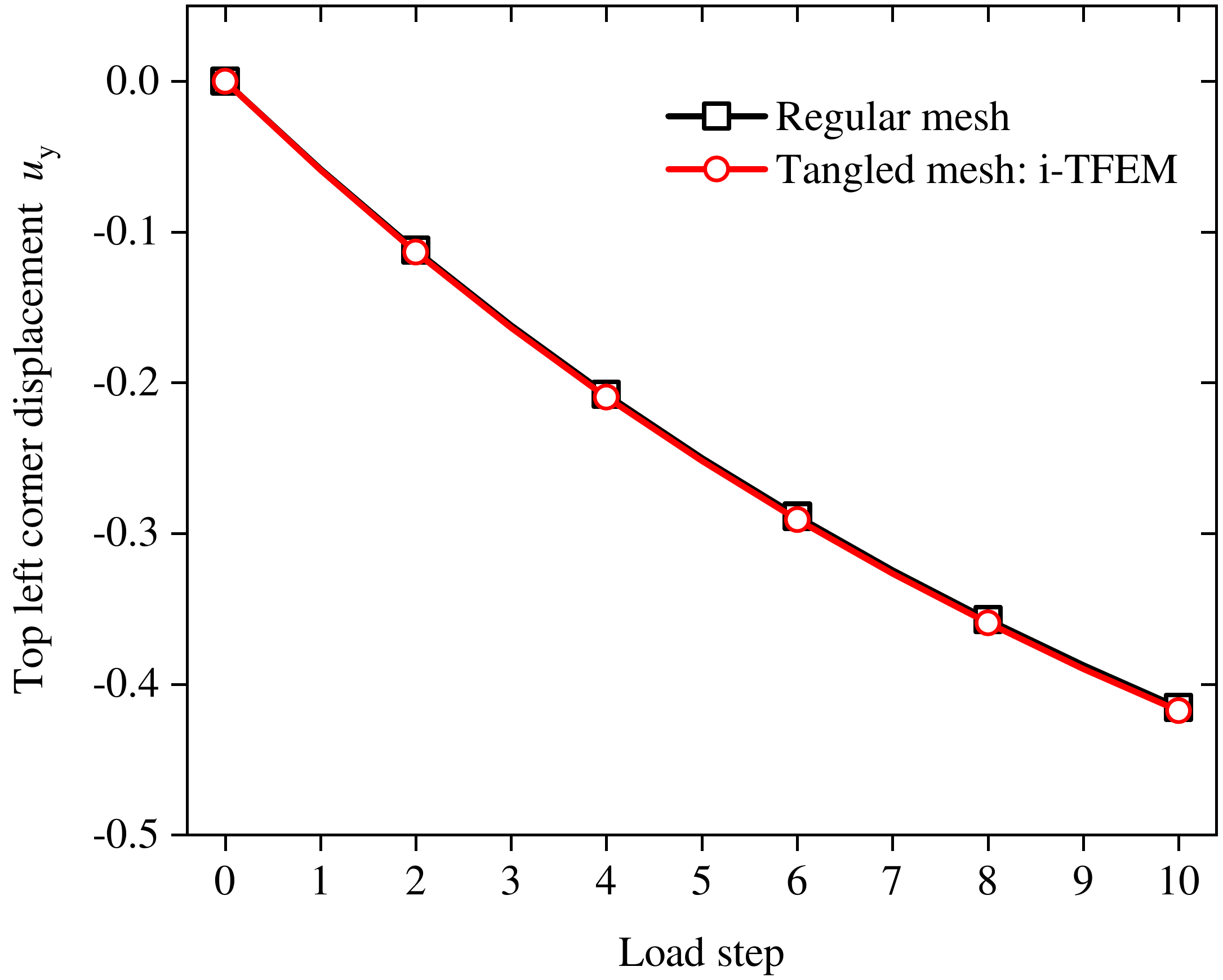}
 	\caption{Vertical displacement versus the load step for the punch problem.}
 	\label{Fig_punchDispVsLoadStep}
 \end{figure}

\begin{figure}[H]
	\begin{subfigure}[c]{.5\textwidth}				
		\centering\includegraphics[width=0.98\linewidth]{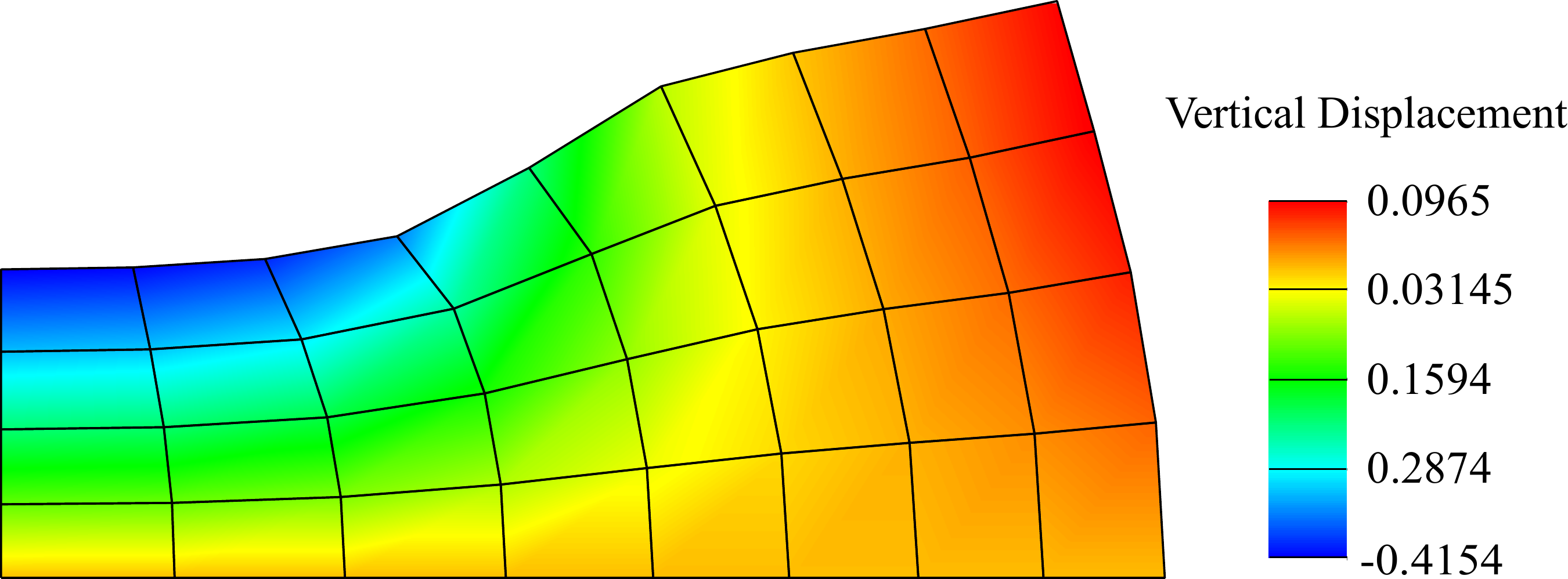}
		\caption{}
		\label{Fig_punchDeformedUntangled}
	\end{subfigure}
	\begin{subfigure}[c]{0.5\textwidth}		
		\centering\includegraphics[width=0.98\linewidth]{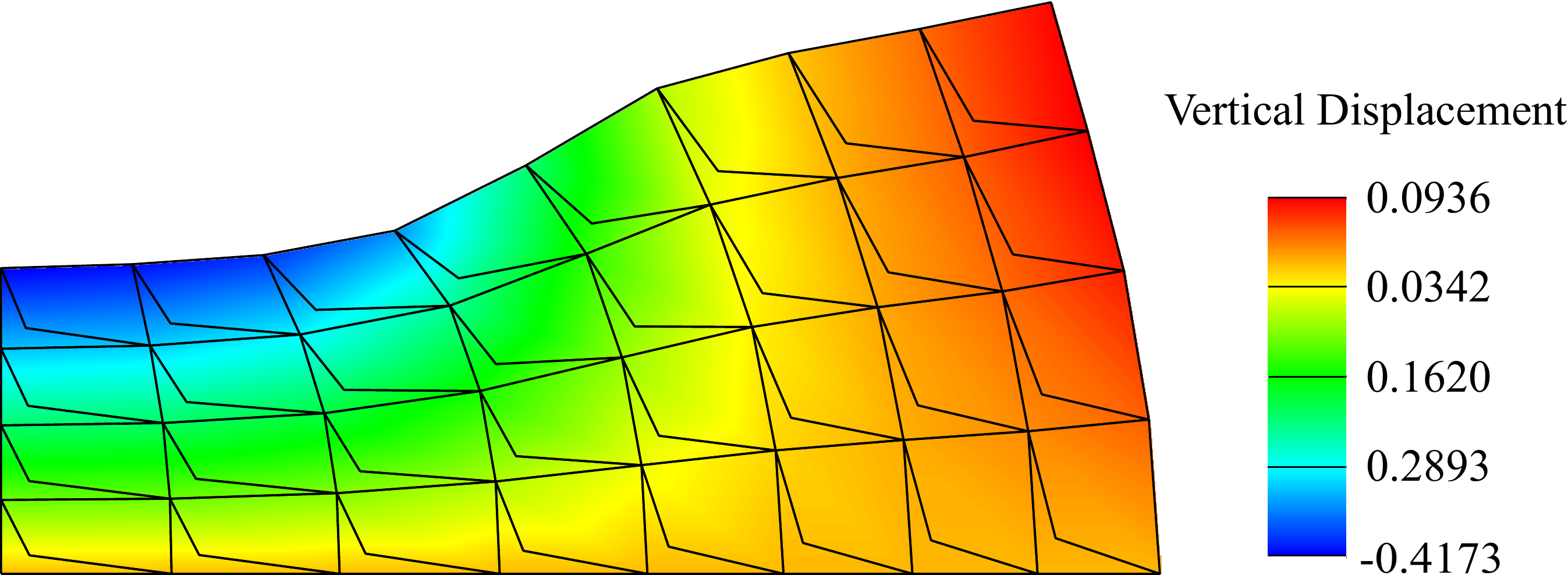}
		\caption{}
		\label{Fig_punchDeformedTangled}
	\end{subfigure}	
	\caption{Deformed configuration for (a) regular mesh via FEM and (b) tangled mesh via i-TFEM for the punch problem.}
\end{figure}

Next, we consider a punch problem \cite{wriggers2017efficient, li2020hyperelastic, de2019serendipity, van2020virtual} with geometric and material nonlinearities. Specifically, compressible isotropic generalized neo-Hookean material model is considered where the strain energy density is given by \cite{zienkiewicz2000finite, bower2009applied}:
\begin{equation}
\Psi_\text{GNH} \left( \boldsymbol{u}   \right)  = \frac{\mu}{2} \left(J_F^{-2/3} tr\boldsymbol{b} -3 \right) + \frac{K}{2}\left(J_F-1 \right)^2 
\label{Eq_neoHookean}
\end{equation}
where $J_F = det \boldsymbol{F} $ and $\boldsymbol{b} = \boldsymbol{F} \boldsymbol{F}^\top $ is the left Cauchy-Green deformation tensor while $\mu = 500$ and  $K = 1700$ are the material parameters (equivalent to shear and bulk moduli respectively in the small strain limit).
A rectangular block is subject to a vertical load $p$ (per unit length) uniformly distributed over the top left half of the block where $p = 1000$ and $H = 1$; see Fig.~\ref{Fig_punchGeometryNMesh} \cite{de2019serendipity}. The top and left sides of the block are fixed in the horizontal direction, while the bottom  is fixed in the vertical direction.  
 Fig.~\ref{Fig_punchDispVsLoadStep} captures the vertical displacement of point A (located at the top left corner) for every load step. The results for the regular mesh and tangled mesh (using i-TFEM) match well. For both meshes, the solution converged  in about 5 Newton iterations for each load step. 
  %The final configurations obtained via i-TFEM for regular and tangled meshes is shown in Fig.~\ref{Fig_punchDeformedUntangled}. 
 Fig.~\ref{Fig_punchDeformedUntangled} and  Fig.~\ref{Fig_punchDeformedTangled}  illustrate the deformed configurations for the regular mesh  and tangled meshes  respectively, after the final load step.  
 
\begin{figure}[H]
	\centering\includegraphics[width=0.5\linewidth]{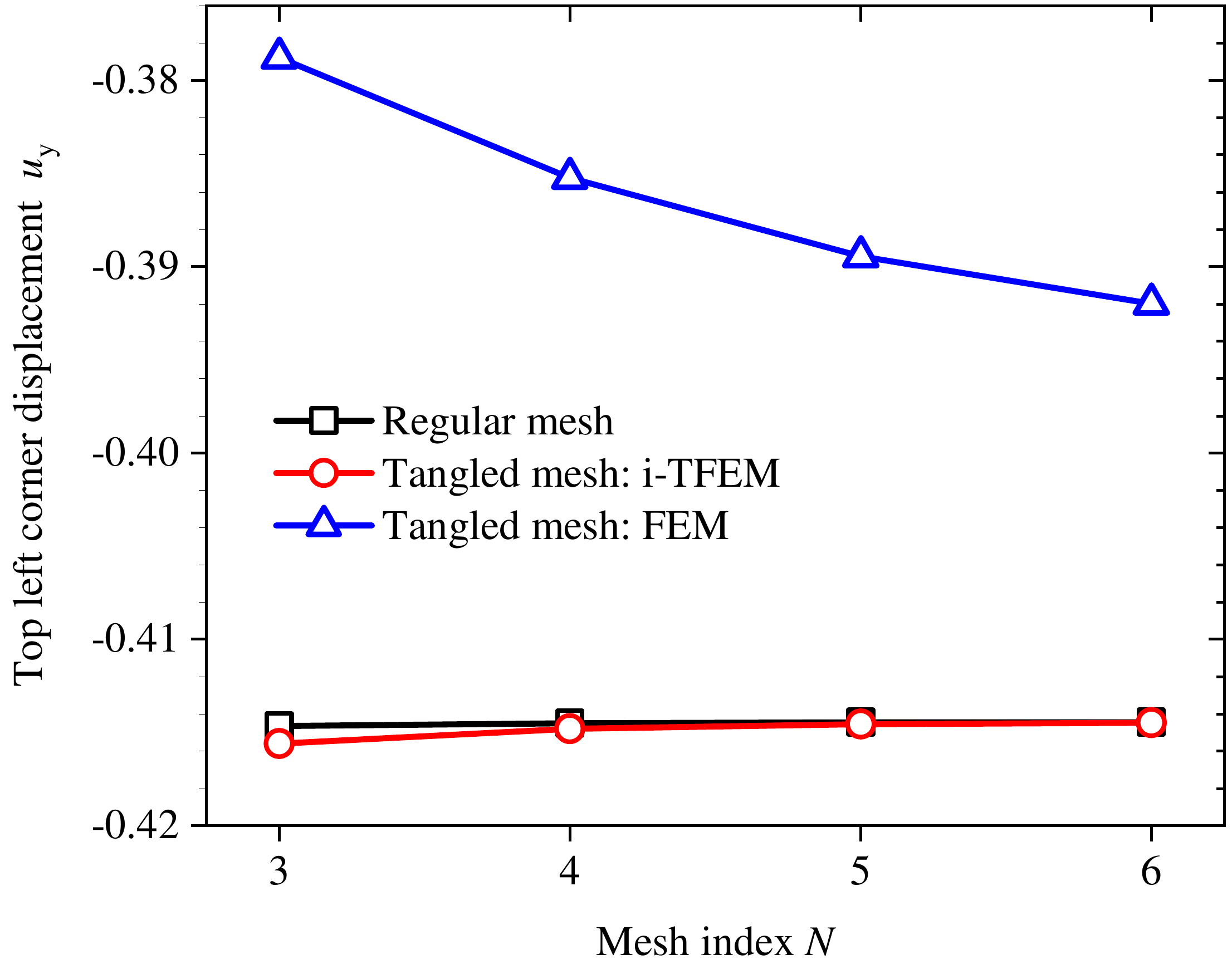}
	\caption{Convergence study for the punch problem.}
	\label{Fig_punchConv}
\end{figure}

\begin{figure}[H]
	\begin{subfigure}[c]{.5\textwidth}				
		\centering\includegraphics[width=0.95\linewidth]{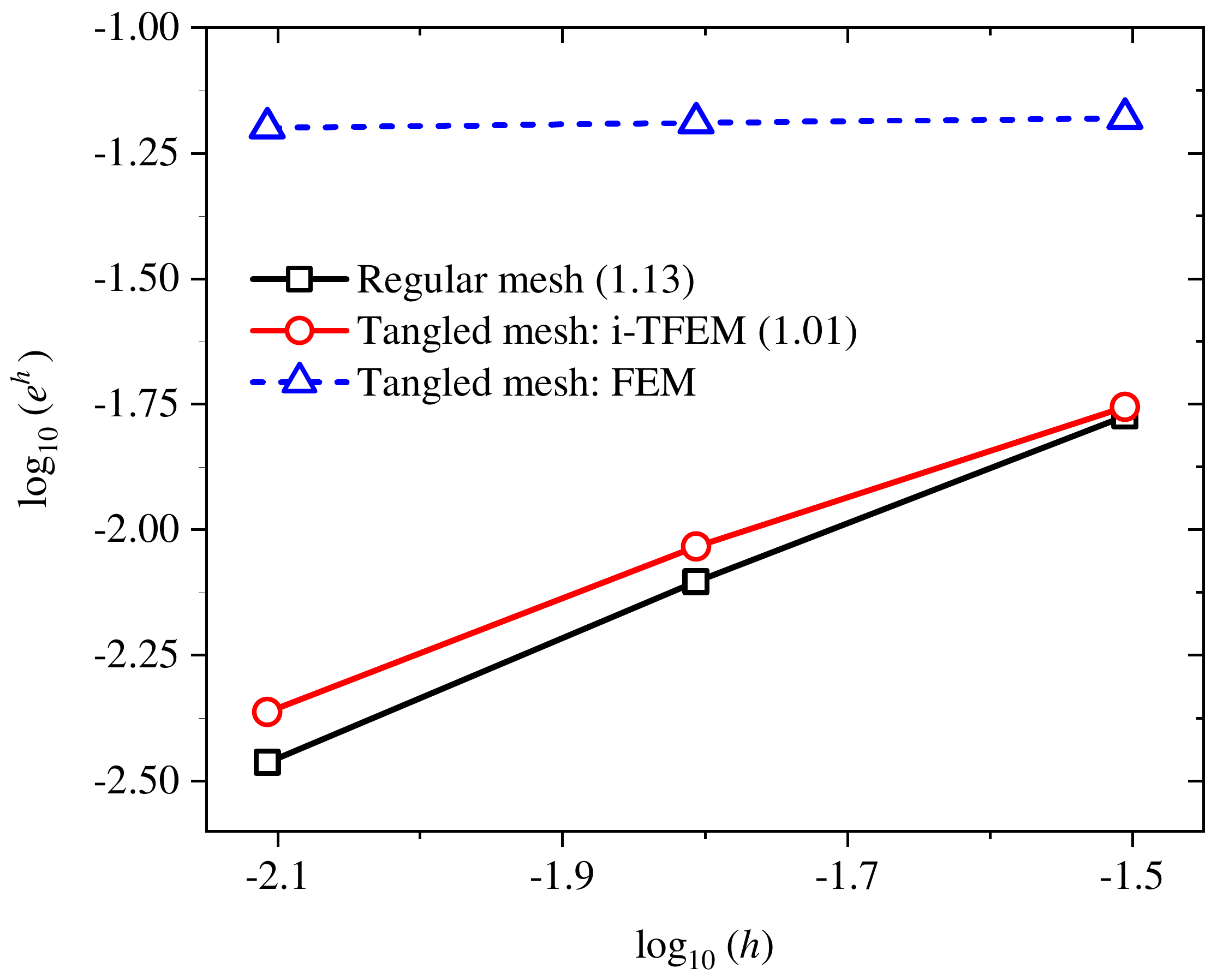}
		\caption{}
		\label{Fig_punchConv_H1}
	\end{subfigure}
	\begin{subfigure}[c]{.5\textwidth}				
		\centering\includegraphics[width=0.95\linewidth]{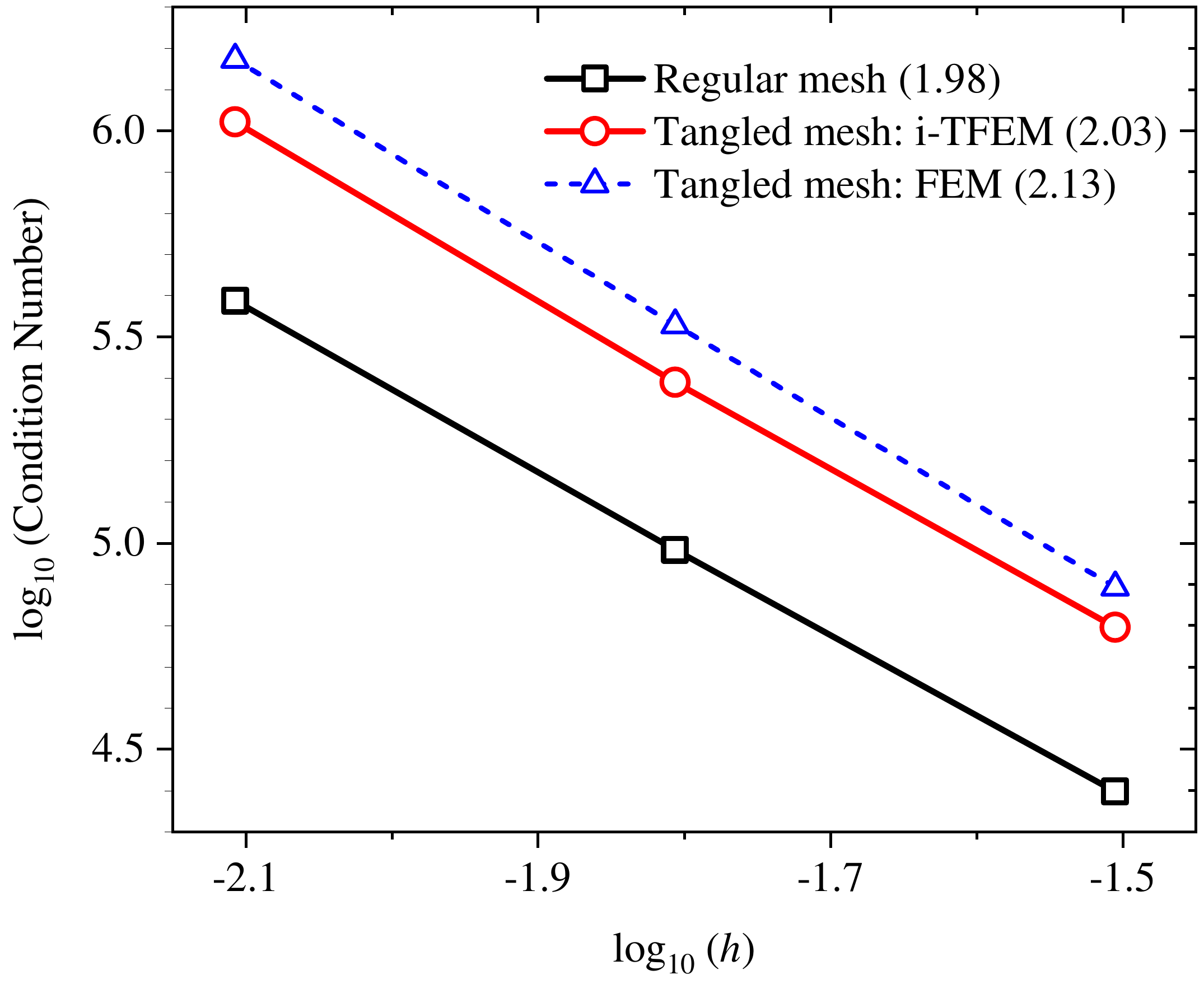}
		\caption{}
		\label{Fig_conditionNumpunch_h}
	\end{subfigure}	
	
	\caption{(a) ${H}^1$ seminorm error versus mesh size and (b) condition number versus mesh size for the punch problem. The convergence rates are provided in the brackets. }
	\label{}	
\end{figure}

To study the convergence, we use the mesh index $N$ where the number of elements in the regular mesh is $ 2^{N+1}\times 2^N$. The regular and tangled meshes with $N = 2$ are shown in Fig.~\ref{Fig_punchGeometryNMesh}.  A convergence study was then carried out as $N$ was varied. The vertical displacement $u_y$ at point A for the two methods is plotted against the mesh index $N$ in Fig.~\ref{Fig_punchConv}. One can observe that the two methods converge to the same solution while standard FEM over the tangled mesh converges to an incorrect solution.
Finally, Fig.~\ref{Fig_punchConv_H1} illustrates the ${H}^1$ seminorm error  over the tangled mesh using FEM and  i-TFEM as well as over the regular mesh. The reference solution is obtained with $N = 8$. Once again, i-TFEM exhibits a convergence rate for the $H^1$ seminorm error (Fig.~\ref{Fig_punchConv_H1}) and the condition number (Fig.~\ref{Fig_conditionNumpunch_h}) similar to that obtained with a non-tangled mesh.		

\subsubsection{Multiple Overlaps}
\begin{figure}[H]
	\centering\includegraphics[width=0.2\linewidth]{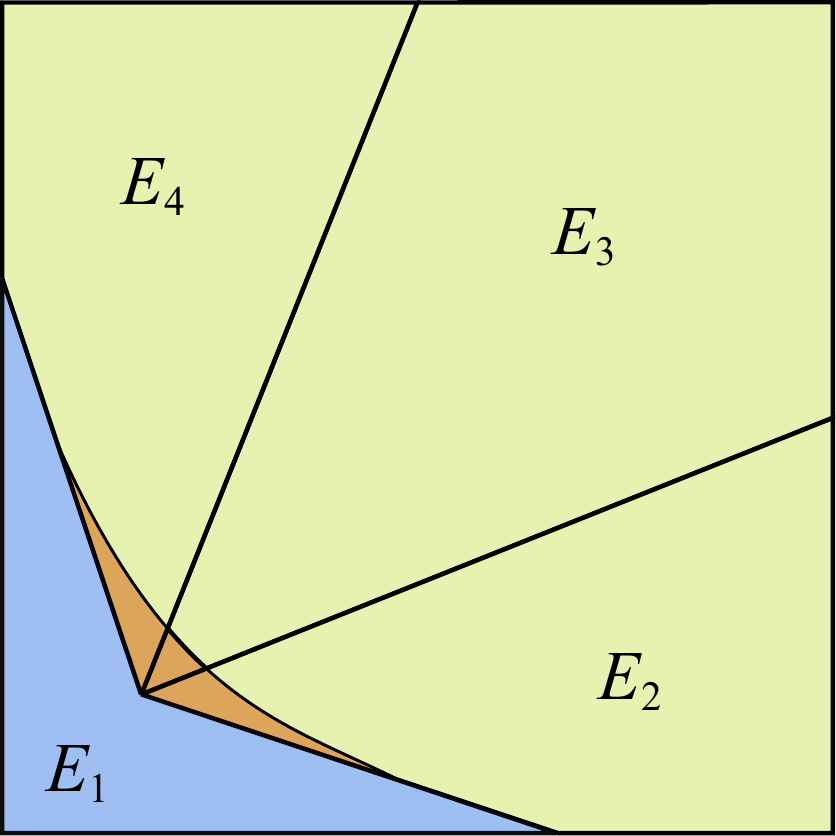}
	\caption{Four-element patch with one element concave.}
	\label{Fig_4elem_fold_colored}
\end{figure}

Thus far, the fold was shared by only one neighboring convex element. However, in practice,  the fold may be shared by multiple convex elements as illustrated in Fig.~\ref{Fig_4elem_fold_colored}. In this case, three convex elements $E_2$, $E_3$ and $E_4$ share the folded region $F_1$. However this does not change the methodology, i.e., the tangent matrices and constraint matrix are computed as before (1) the tangent matrix $\boldsymbol{K}_{convex}^{t}$ and residual vector $\boldsymbol{R}^u_{convex}$  corresponding to the convex elements   are computed using the three convex elements,  (2) while $\widehat{\boldsymbol{K}}_{concave}^{t}$ and $\widehat{\boldsymbol{R}}^u_{concave}$ are computed using the parametric space associated with $ \widehat{E}_1$,  and (3) the constraint matrix is computed using the entire fold $F_1$.

\begin{figure}[H]
	\begin{subfigure}[c]{.45\textwidth}				
		\centering\includegraphics[width=0.9\linewidth]{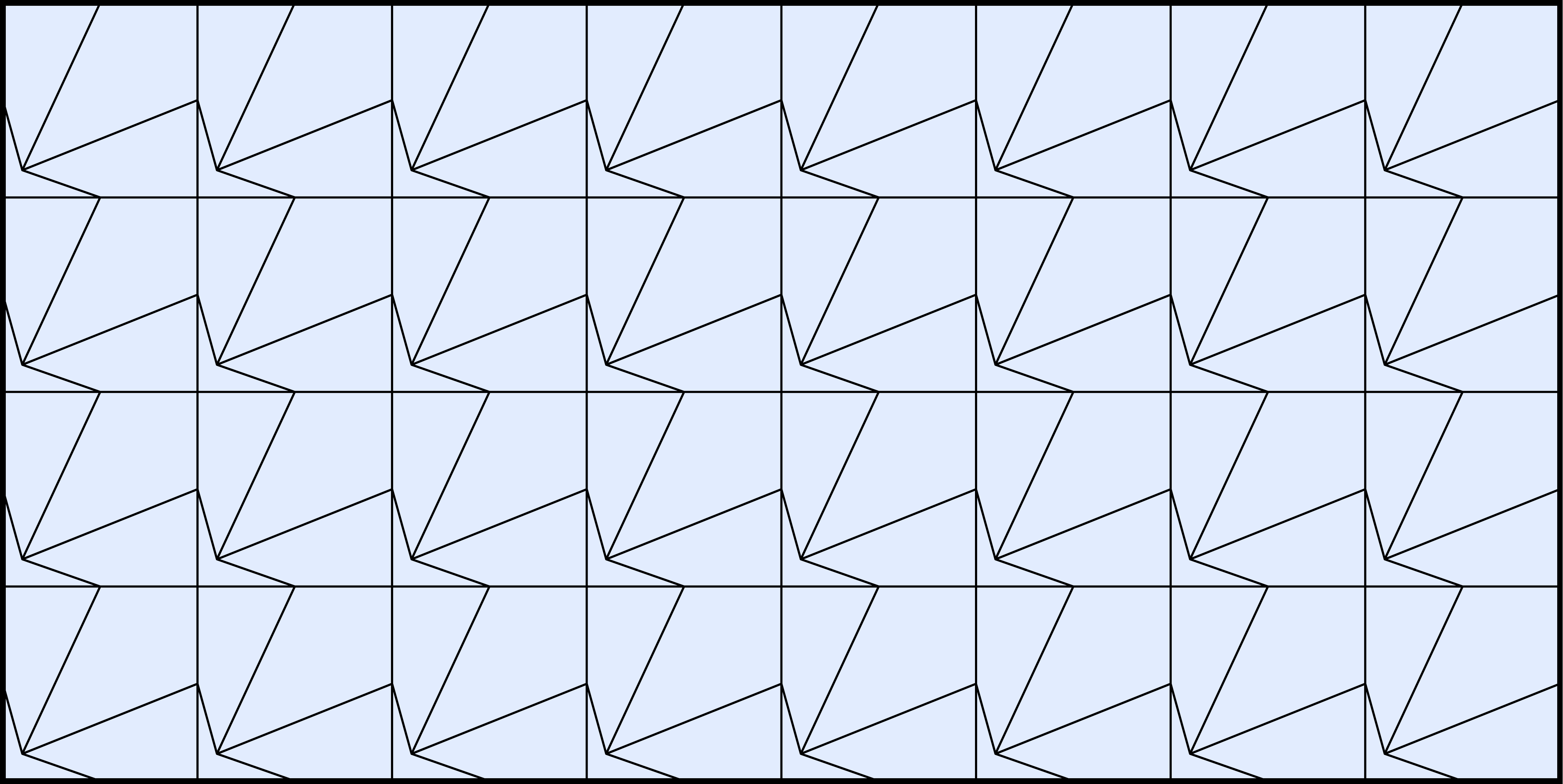}
		\caption{}
		\label{Fig_punch_VDIsp_multipleOverlap_TFEM3x3}
	\end{subfigure}
	\begin{subfigure}[c]{.55\textwidth}				
		\centering\includegraphics[width=1\linewidth]{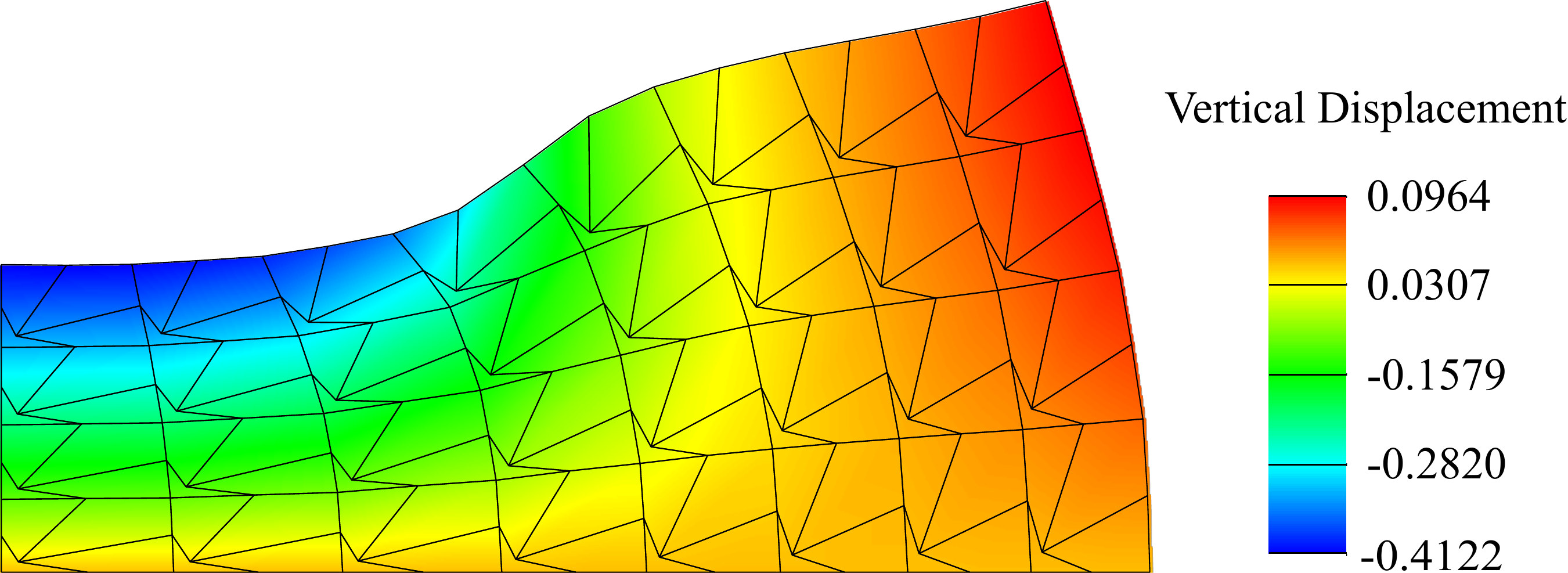}
		\caption{}
		\label{Fig_punch_4elelmVdisp_deformed}
	\end{subfigure}	
	
	\caption{(a) $N = 3$  tangled mesh  (b) Final deformed shape obtained via i-TFEM. }
	\label{}	
\end{figure}

\begin{figure}[H]
	\begin{subfigure}[c]{.5\textwidth}				
		\centering\includegraphics[width=0.95\linewidth]{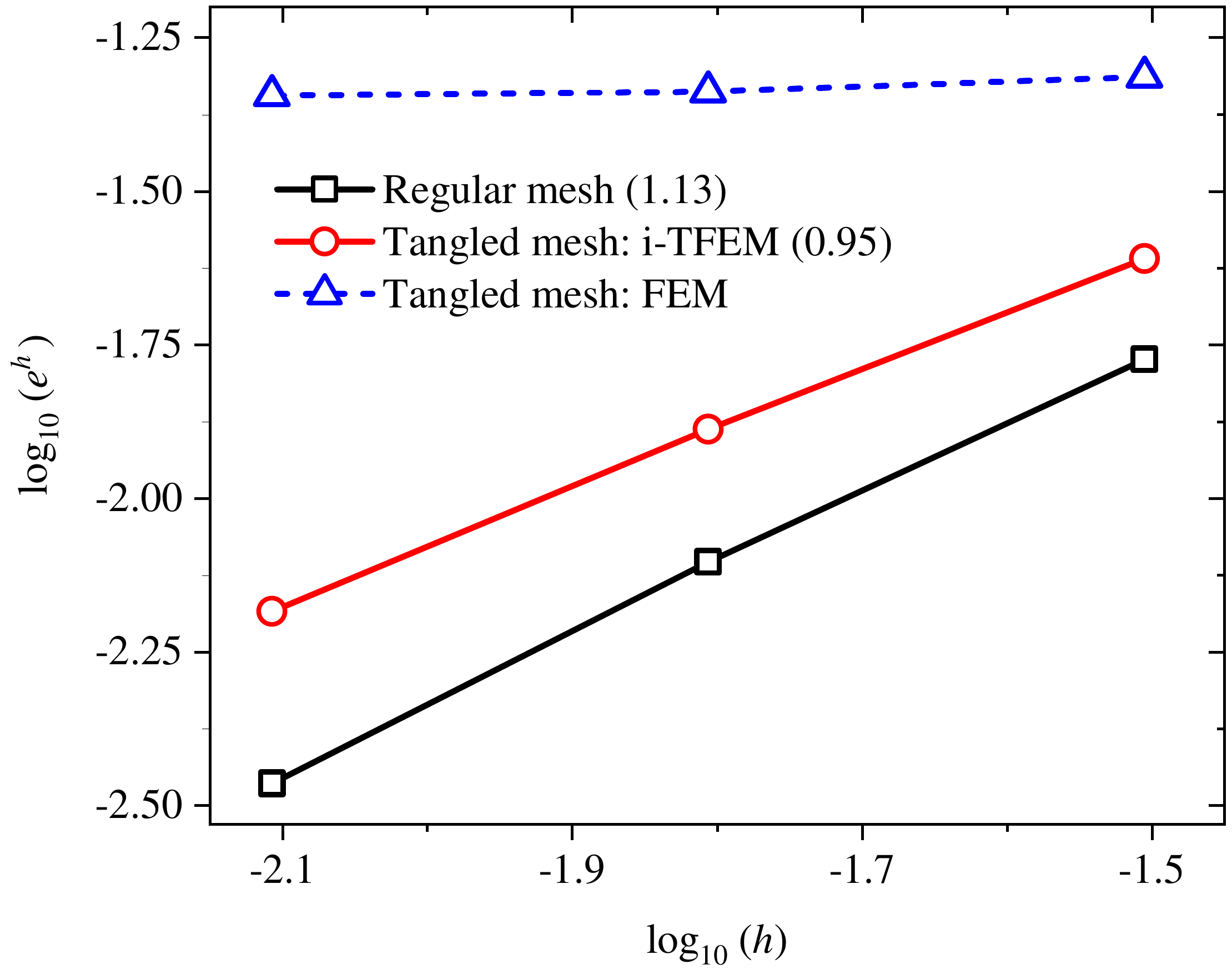}
		\caption{}
		\label{Fig_punch_4elelmGr_H1}
	\end{subfigure}
	\begin{subfigure}[c]{.5\textwidth}				
		\centering\includegraphics[width=0.95\linewidth]{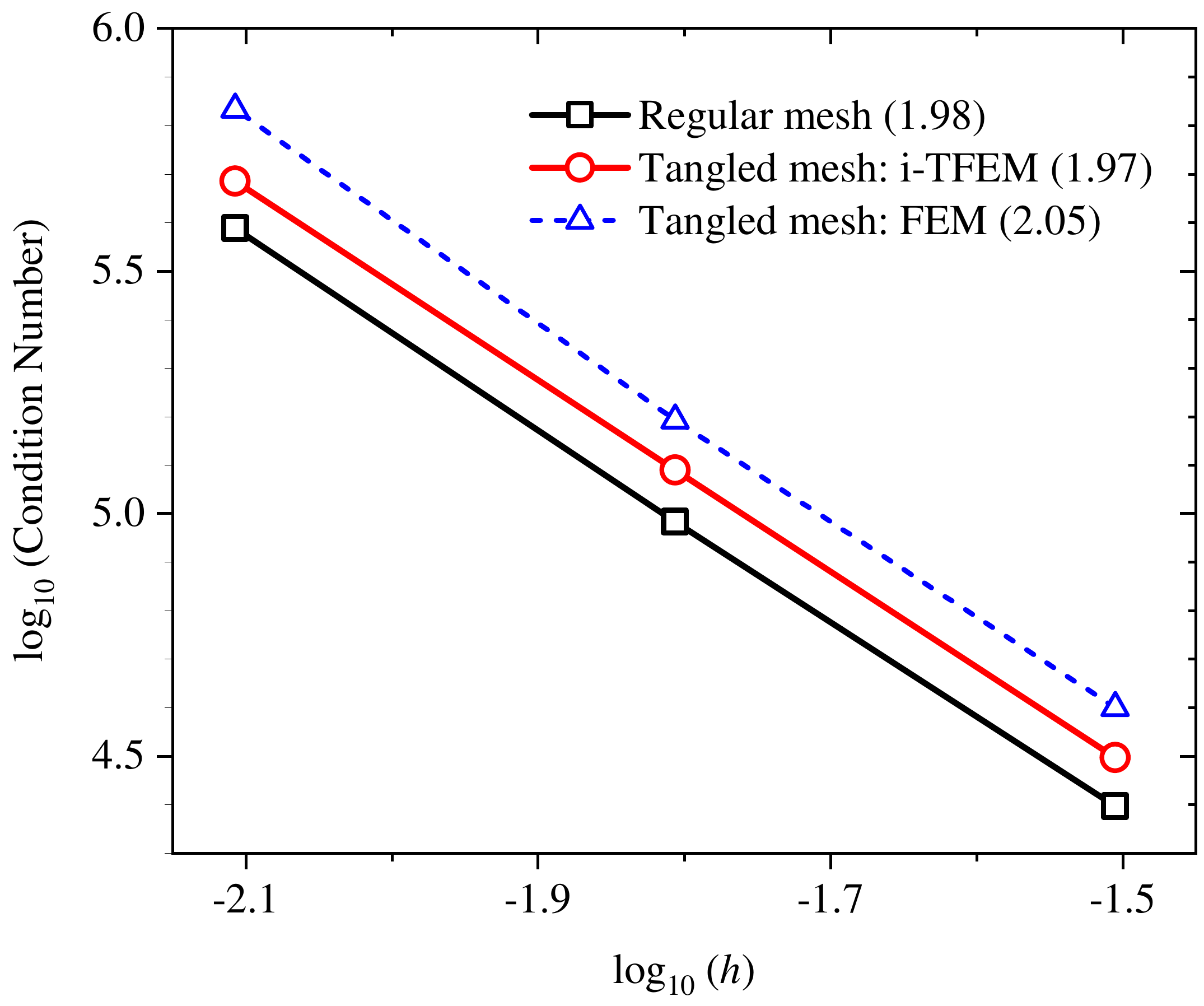}
		\caption{}
		\label{Fig_punch_4elelmGr_CN}
	\end{subfigure}	
	
	\caption{(a) ${H}^1$ seminorm error versus mesh size and (b) condition number versus mesh size for the punch problem with the four-element patch as the repeating unit. The convergence rates are provided in the brackets.}
	\label{}	
\end{figure}

Here, we consider a tangled mesh (see Fig.~\ref{Fig_punch_VDIsp_multipleOverlap_TFEM3x3}) where the basic repeating unit is the patch shown in Fig.~\ref{Fig_4elem_fold_colored}. The problem described in Fig.~\ref{Fig_punchGeometryNMesh} is used as a case study. The final deformed configuration obtained via i-TFEM is shown in Fig.~\ref{Fig_punch_4elelmVdisp_deformed}.
Next, to study the convergence, we use the mesh index $N$ where the number of elements is $ 2^{N+1}\times 2^N$. 
\begin{comment}
The vertical displacement $u_y$ at point A is plotted against the mesh index $N$ in Fig.~\ref{Fig_punchConv} with i-TFEM and FEM over the tangled mesh as well as the regular mesh. One can observe that i-TFEM converges to the expected solution as opposed to the standard FEM.

\begin{figure}[H]
	\centering\includegraphics[width=0.5\linewidth]{punch_Vdisp_conv_4elem.pdf}
	\caption{Convergence study for the punch problem with the four-element patch as the repeating unit.}
	\label{Fig_punch_4elelmGr_conv}
\end{figure}
\end{comment}
Fig.~\ref{Fig_punch_4elelmGr_H1} illustrates the ${H}^1$ seminorm error  over the tangled mesh using FEM and  i-TFEM as well as over the regular mesh. The reference solution is obtained with $N = 8$. Once again, i-TFEM exhibits an optimal convergence rate (Fig.~\ref{Fig_punch_4elelmGr_H1}). Moreover, the condition number for i-TFEM is comparable to FEM (Fig.~\ref{Fig_punch_4elelmGr_CN}).		

\subsubsection{Incompressibility}

Next, we evaluate the performance of i-TFEM under near-incompressibility. Here, the F-bar method \cite{neto2005f} is employed  to prevent locking, both in FEM and i-TFEM. 
The punch problem (Fig.~\ref{Fig_punchGeometryNMesh}) is considered with the neo-Hookean material parameters  $\mu = 500$ and $K = 5 \times 10^6$ which results in Poisson's ratio $\nu = 0.49995$ (near-incompressibility). The problem is solved over the regular (tangle-free) mesh (shown in Fig.~\ref{Fig_punchGeometryNMesh}) and the tangled mesh (Fig.~\ref{Fig_punch_VDIsp_multipleOverlap_TFEM3x3}). The size of the mesh is governed by the index $N$ as discussed previously. 
Fig.~\ref{Fig_incomp_punch_N6}  shows the deformed configuration obtained over the regular mesh and via i-TFEM over tangled mesh with the mesh index $N = 4$. Observe the close agreement between the two solutions.  

\begin{figure}[H]
	\begin{subfigure}[c]{.5\textwidth}				
		\centering\includegraphics[width=0.98\linewidth]{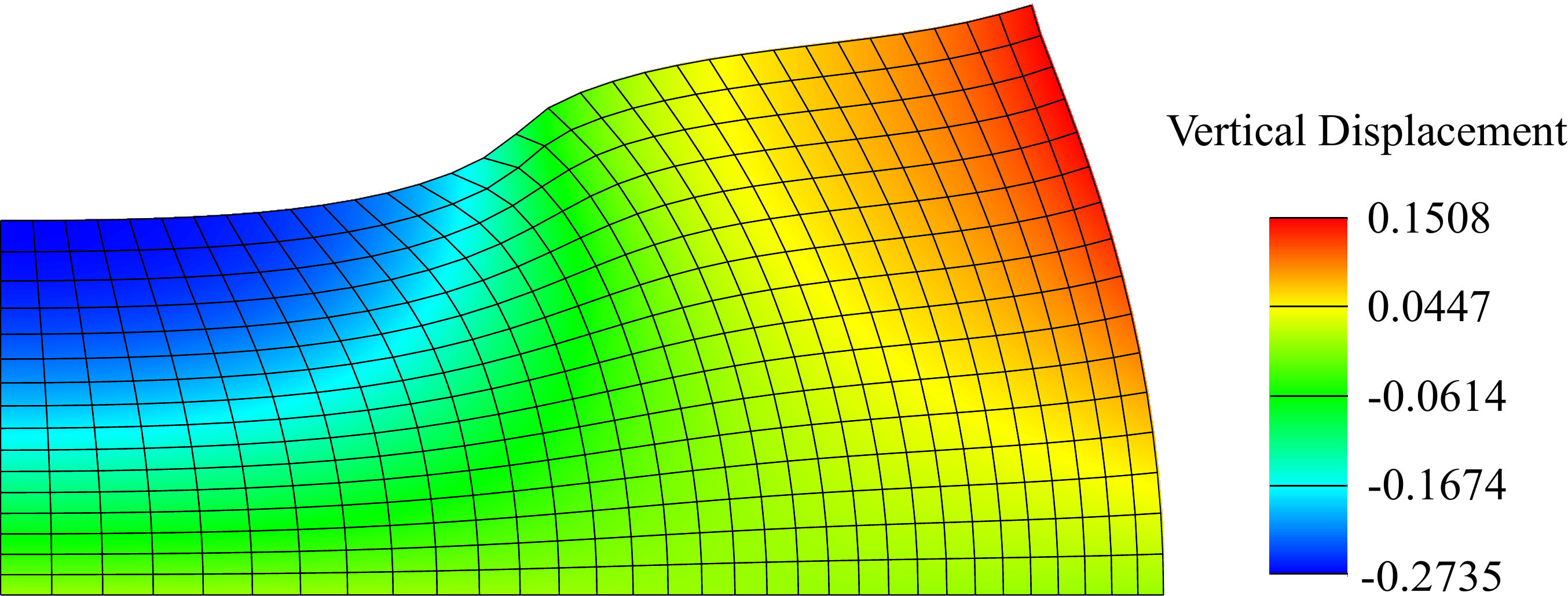}
		\caption{}
		\label{Fig_punchDeformedUntangled_incomp}
	\end{subfigure}
	\begin{subfigure}[c]{0.5\textwidth}		
		\centering\includegraphics[width=0.98\linewidth]{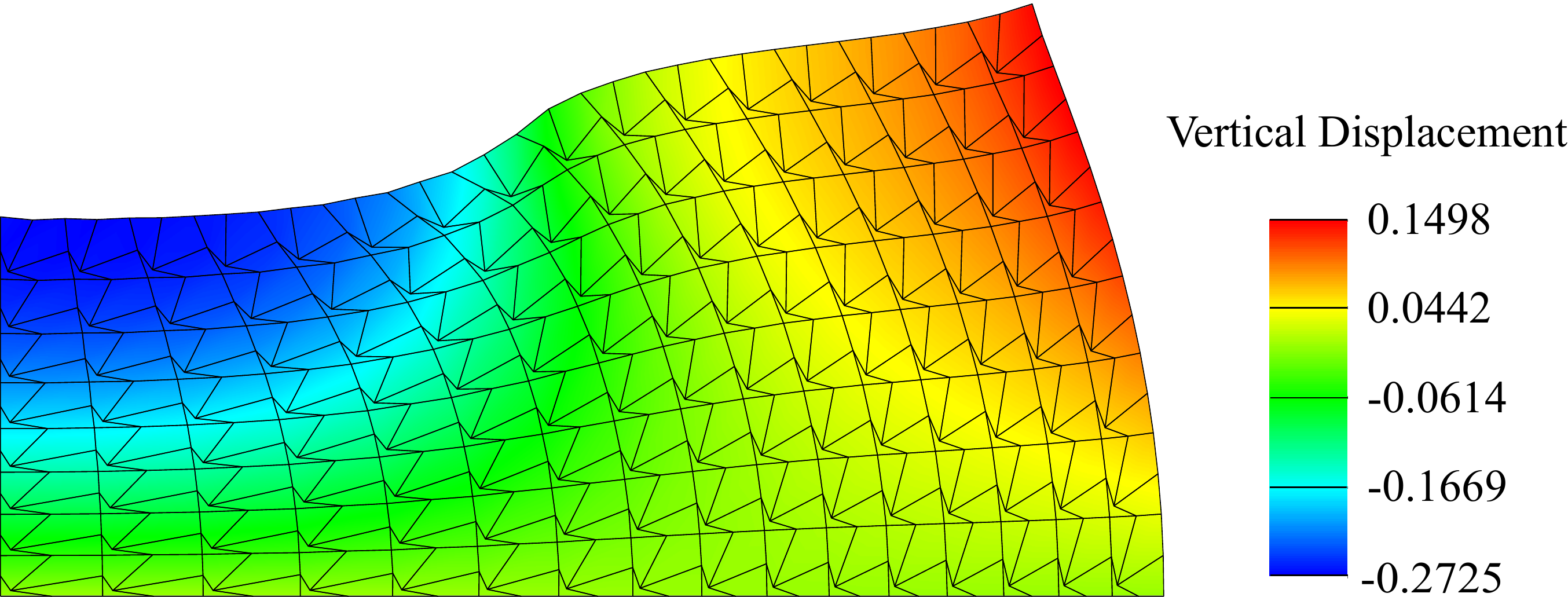}
		\caption{}
		\label{Fig_punchDeformedTangled_incomp}
	\end{subfigure}	
	\caption{Deformed configuration for (a) regular mesh via FEM and (b) tangled mesh via i-TFEM for the punch problem.}
\end{figure}
\begin{figure}[H]
	\centering\includegraphics[width=0.5\linewidth]{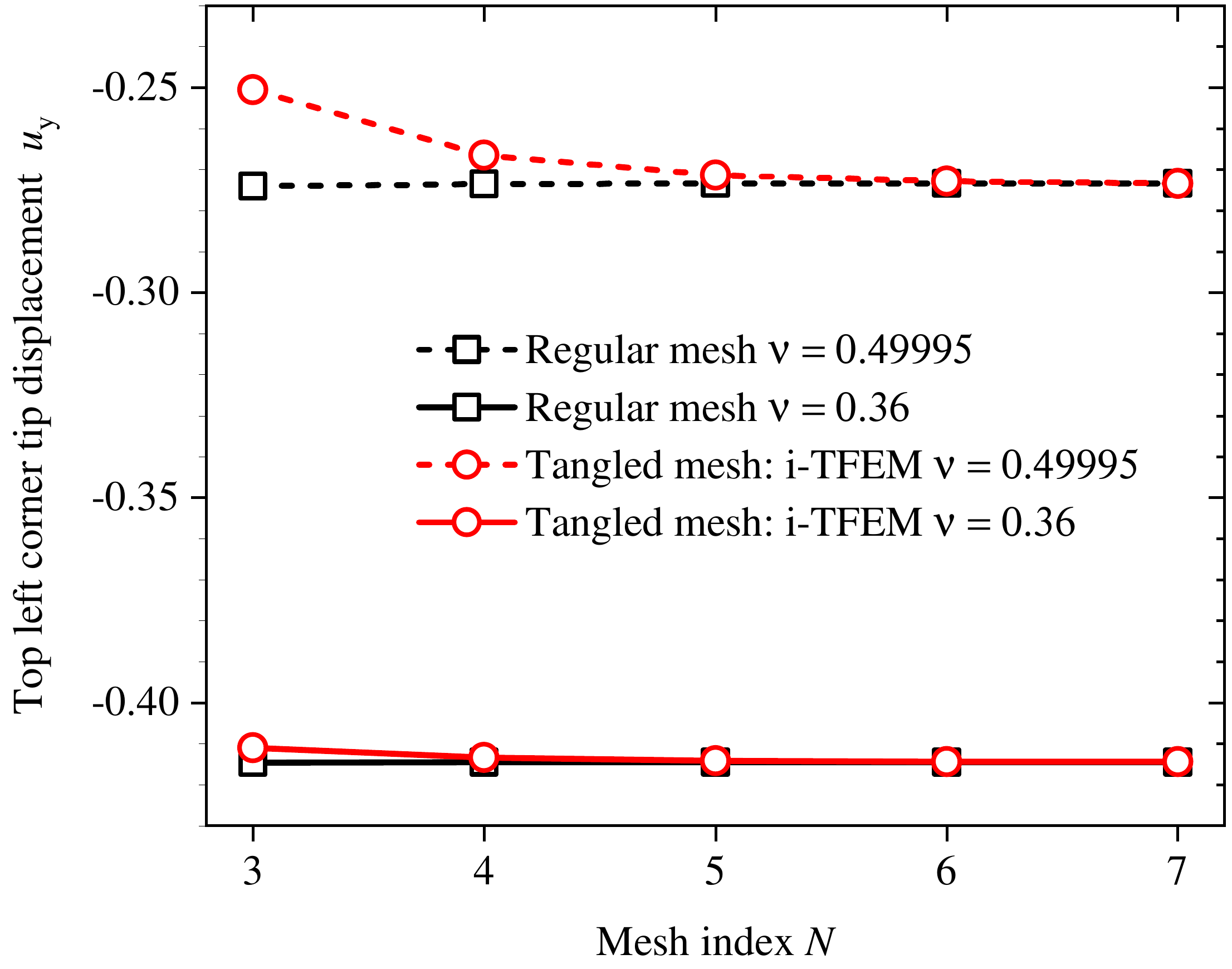}
	\caption{Convergence study for the punch problem with the four-element patch as the repeating unit.}
	\label{Fig_punch_4elelmGr_conv_2PoissonRatio}
\end{figure}

\begin{figure}[H]
	\centering\includegraphics[width=0.75\linewidth]{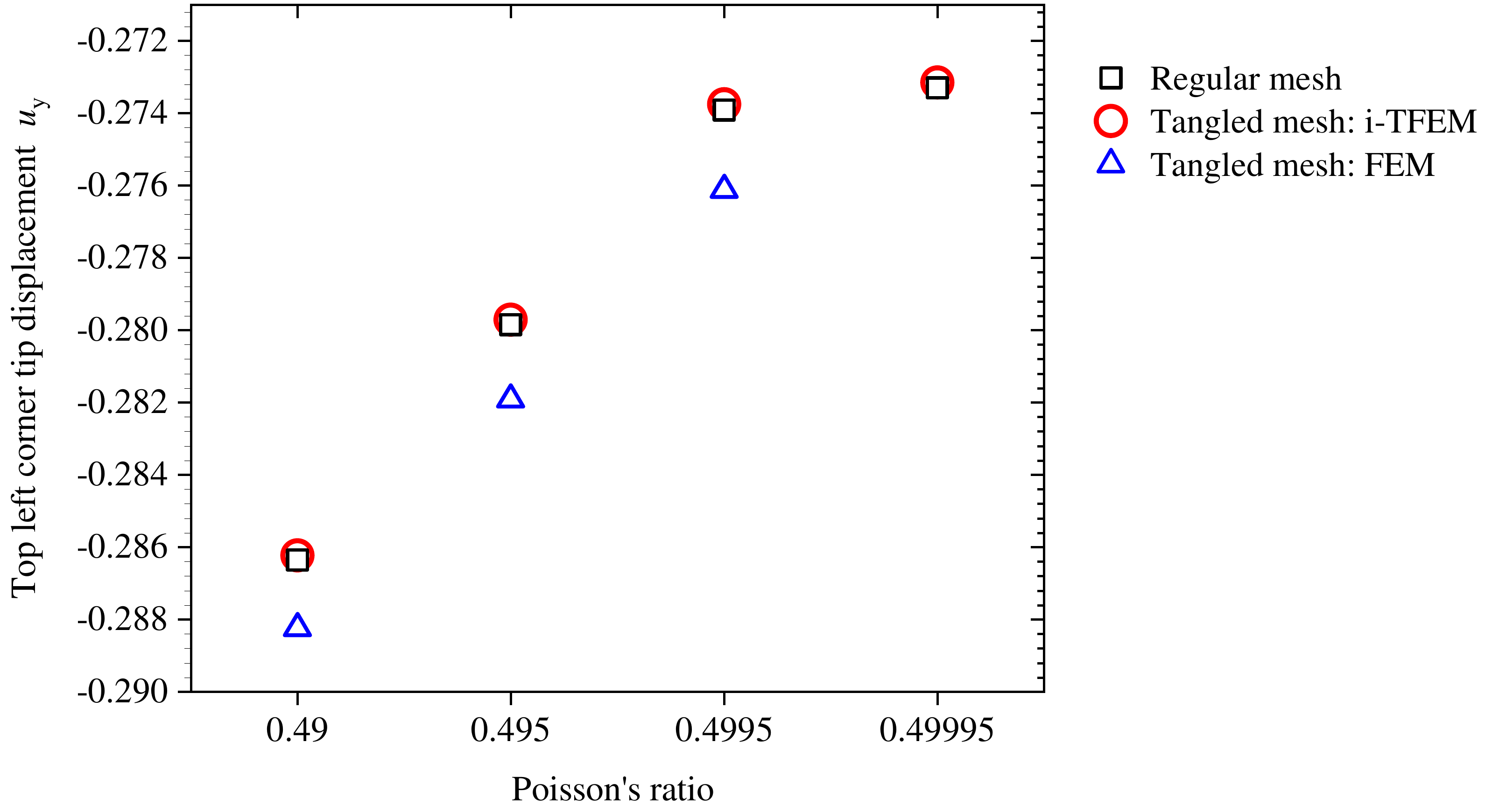}
	\caption{Tip vertical displacement versus different Poisson ratios for the punch problem.}
	\label{Fig_incomp_punch_N6}
\end{figure}

	
Next, we study the convergence of i-TFEM for two Poisson's ratio values: (a) 0.49995 and (b) 0.36.  Note that for Poisson's ratio  = 0.36, the F-bar method is not required. The vertical displacement $u_y$ at point A is plotted against the mesh index $N$ in Fig.~\ref{Fig_punch_4elelmGr_conv_2PoissonRatio} for the two methods. One can observe that i-TFEM converges to same solution as FEM in both cases. 

Finally, we compare the performance of FEM and i-TFEM for different Poisson's ratios close to 0.5. For this experiement, we consider a regular mesh and a tangled mesh with $N = 7$.  Fig.~\ref{Fig_incomp_punch_N6} plots the vertical displacement at point A obtained by FEM and i-TFEM over tangled and untangled meshes for four values of Poisson's ratio: 0.49, 0.495, 0.4995 and 0.49995. As opposed to FEM, the solutions obtained via i-TFEM over the tangled mesh converge to (approximately) the same value as that obtained over the untangled mesh.  Note that the FEM solution  over the tangled mesh at $\nu = 0.49995$ is not provided since the Newton-Raphson method did not converge even with a sufficiently large number of load steps. 

\subsection{Thin Beam}
\begin{figure}[H]
	\begin{subfigure}[c]{1\textwidth}				
		\centering\includegraphics[width=1\linewidth]{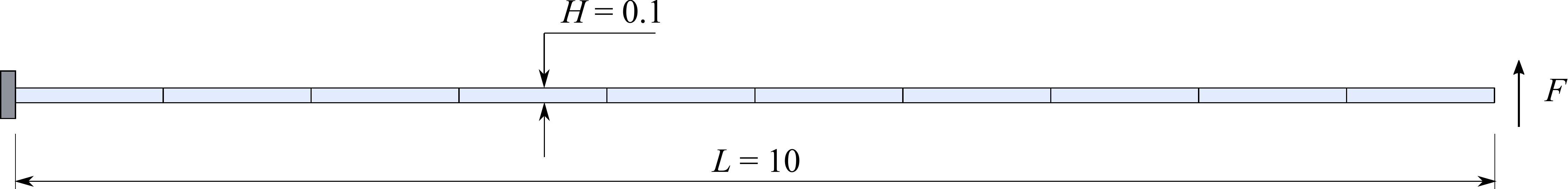}
		\caption{}
		\label{Fig_thinBeamGeom}
	\end{subfigure}
	\begin{subfigure}[c]{1\textwidth}		
		\centering\includegraphics[width=1\linewidth]{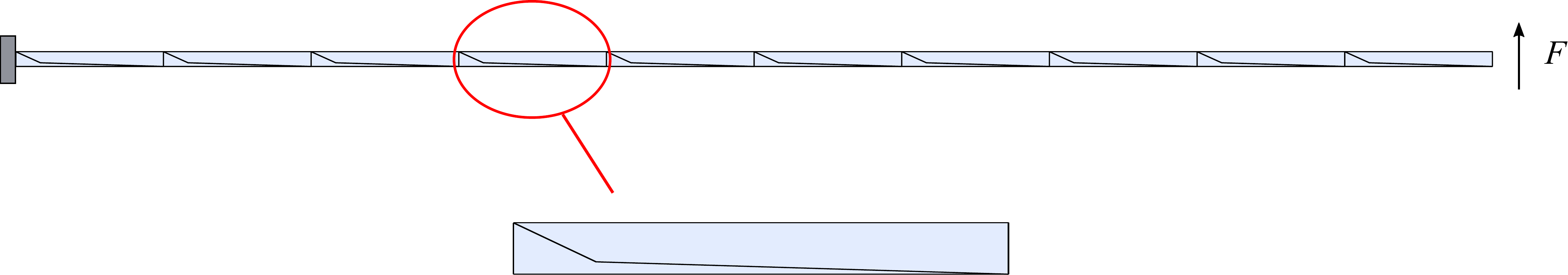}
		\caption{}
		\label{Fig_thinBeamMeshTangled}
	\end{subfigure}	
	\caption{ (a) Thin beam geometry and boundary conditions with  regular mesh  and (b) the corresponding tangled mesh. The repeating unit for the tangled mesh is zoomed in.}
\end{figure}

\begin{figure}[H]
	\centering\includegraphics[width=1\linewidth]{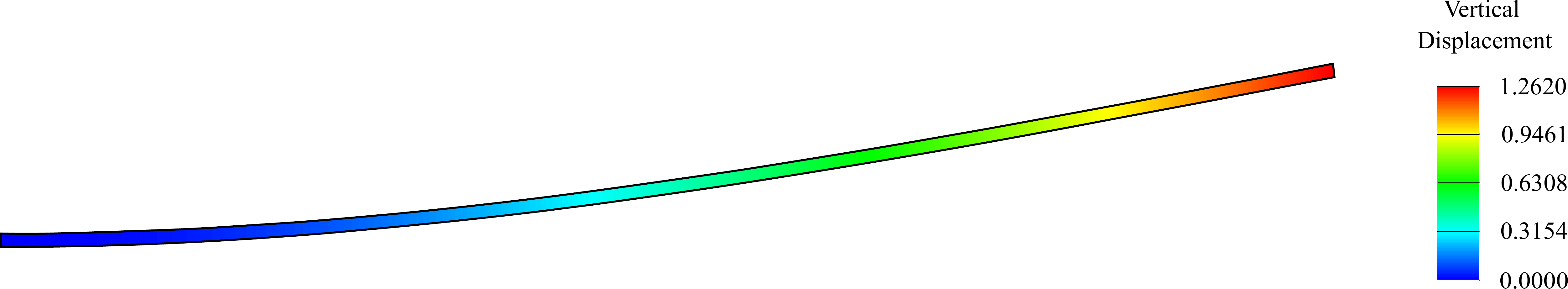}
	\caption{Final deformed shape using i-TFEM with $N = 3$ tangled mesh.}
	\label{Fig_thinBeam_VDisp_TFEM3}
\end{figure}

In this example, a beam undergoing large deflections is considered  \cite{li2020hyperelastic, wriggers2017efficient, de2019serendipity}. Specifically, a beam (see Fig.~\ref{Fig_thinBeamGeom}) with a length-to-height ratio $L/H = 100$ is fixed at the left end and subjected to a vertical load $F = 0.1$ at the right end. The material parameters of the Neo-Hookean  model (Eq.~\ref{Eq_neoHookean}) are $K = 16000$ and $\mu = 6000$. The regular and tangled meshes are shown in Fig.~\ref{Fig_thinBeamGeom} and  Fig.~\ref{Fig_thinBeamMeshTangled}, respectively.
The number of elements in the mesh is governed by the mesh index $N$. For the regular mesh, the number of elements in the horizontal direction is  given as $(10\times2^N)$  while in the vertical direction, the number of elements is given by $2^N$  The regular mesh in Fig 23 a corresponds to $N = 0$. To obtain the corresponding tangled mesh, each element is divided into a concave and a convex element. Hence the total number of elements in the tangled mesh is $2 \times (10\times2^N) \times 2^N$. The final deformed configuration of the beam obtained via i-TFEM over the tangled mesh with $N = 3$ is shown in Fig.~\ref{Fig_thinBeam_VDisp_TFEM3}.

\begin{figure}[H]
	\centering\includegraphics[width=0.5\linewidth]{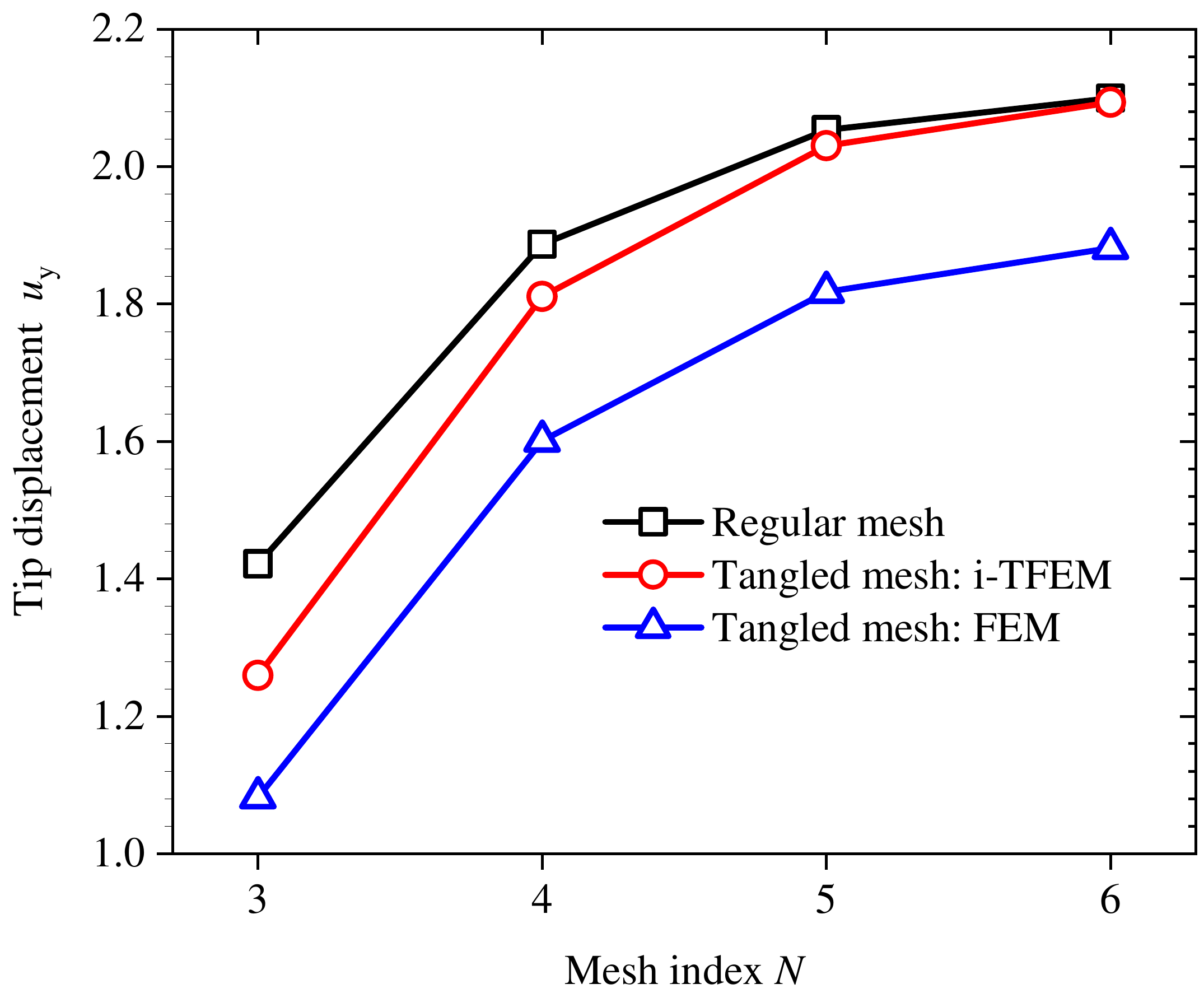}
	\caption{Convergence study for thin beam bending problem.}
	\label{Fig_ThinBeamConv}
\end{figure}

To study the convergence, the vertical displacement at the top right corner of the beam is considered.   Fig.~\ref{Fig_ThinBeamConv} plots the convergence of the regular mesh as well as the tangled mesh with FEM and i-TFEM. It can be seen that with i-TFEM, the solution converges to the same value as that obtained using the regular mesh.

\subsection{Application: Aircraft model}
\label{Sec_AircraftModel}
\begin{figure}[H]
	\setlength{\figSize}{0.5\linewidth}
	\centering			
	\begin{subfigure}[b]{0.9\figSize}
		\centering
		\begin{tikzpicture}[spy using outlines={red, circle, magnification=3, size=1.2cm,
				connect spies}]			
			
			\node{\includegraphics[width=\linewidth]{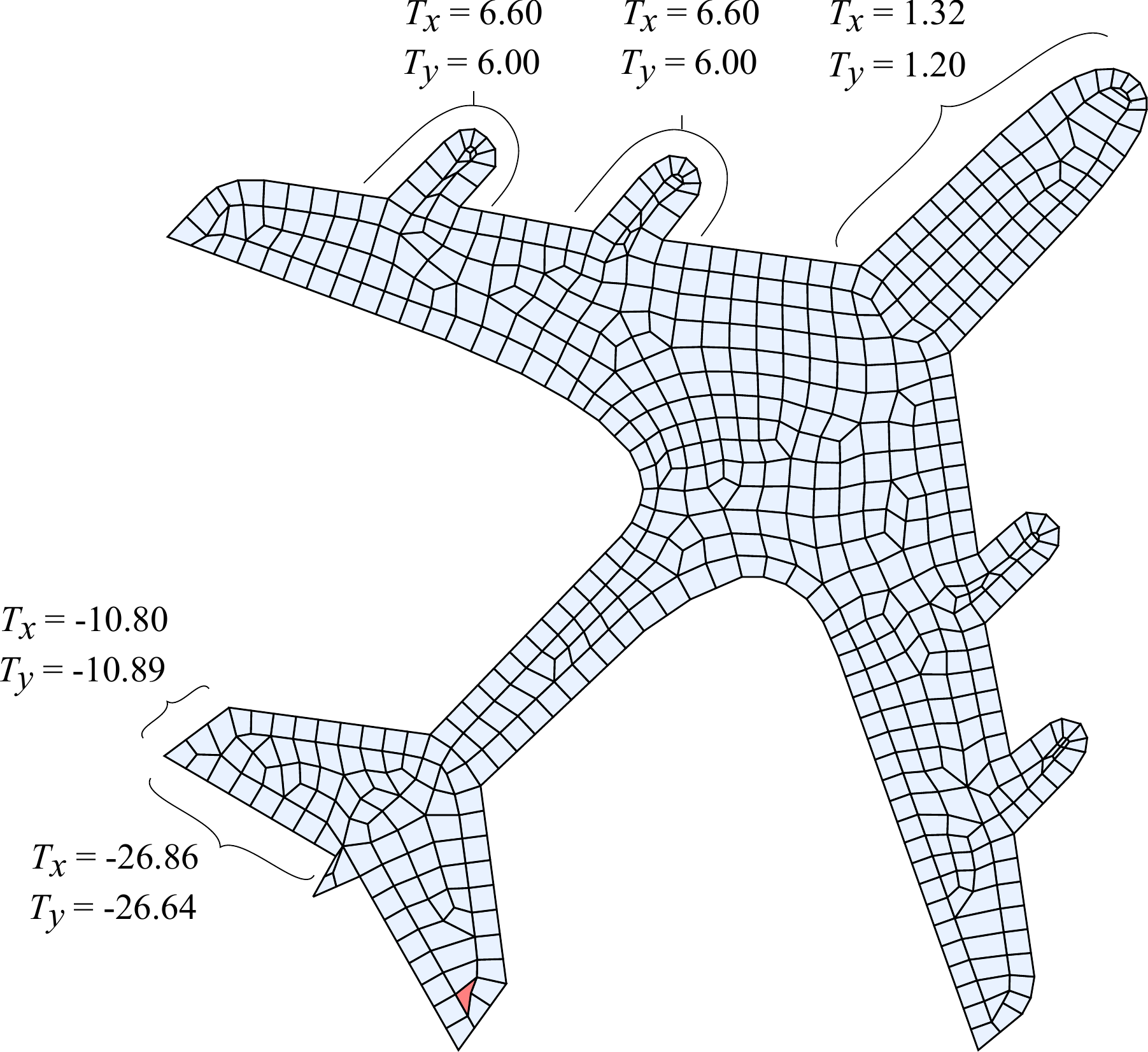}};			
			
			\spy on (-0.7,-3.0) in node [left] at (-1.2,0.6);
		\end{tikzpicture}
		\caption{}
		\label{Fig_aircraft_mesh}
	\end{subfigure}
	\begin{subfigure}[b]{\figSize}
		\centering\includegraphics[width=1\linewidth]{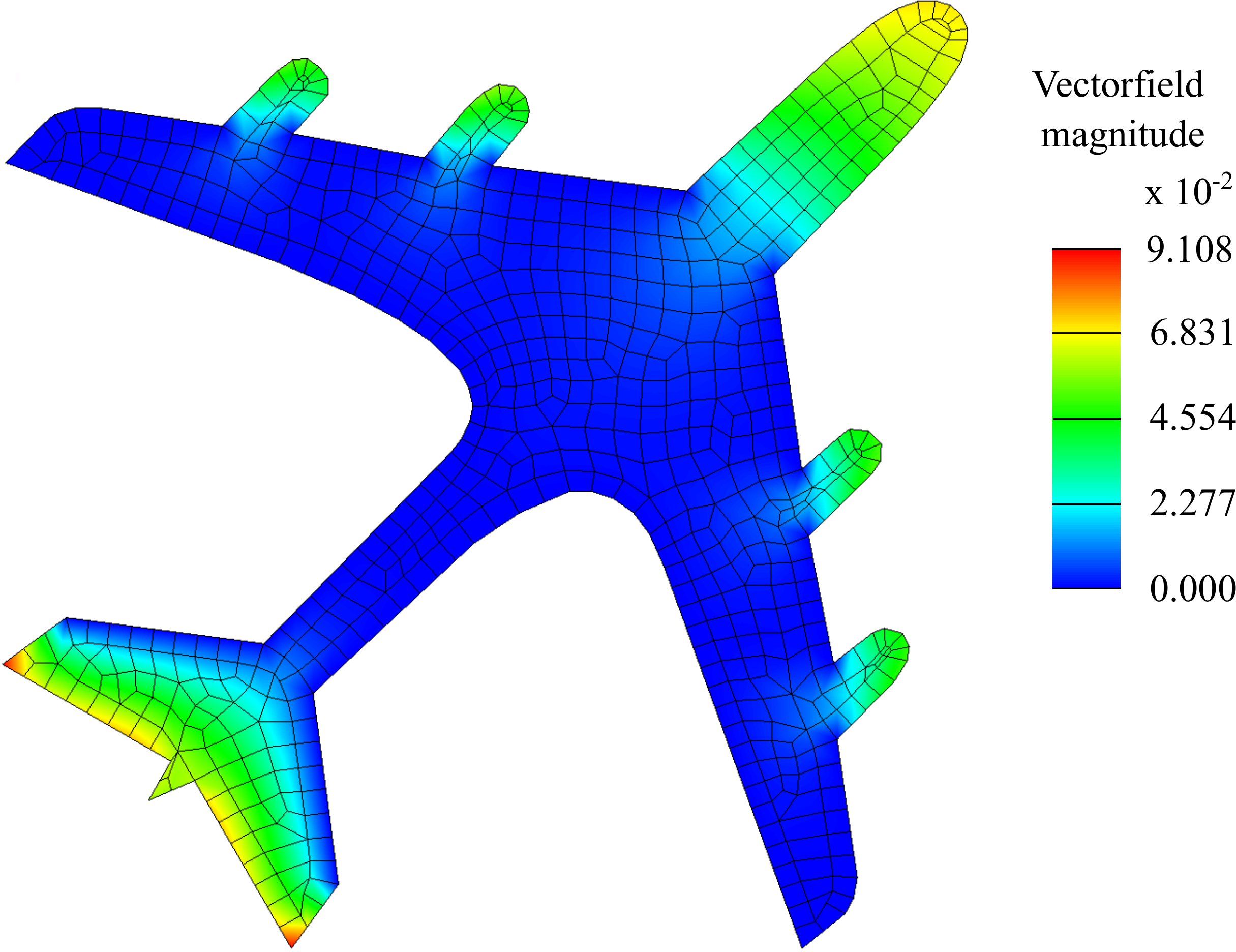}
		\caption{}
		\label{Fig_aircraft_pp}
	\end{subfigure}
	\caption{(a) Mesh for an aircraft model, with one concave element. (b)  i-TFEM solution.}
	\label{Fig_aircraftMeshzoom}
\end{figure}

An example of a practical tangling problem is shown in Figure~\ref{Fig_aircraft_mesh}, where the quadrilateral mesh for an aircraft model was created using the quad mesher proposed in \cite{sarrate2000efficient}. One quad element (out of 600) was found to be concave, for this particular model. The material parameters of the St. Venant-Kirchhoff model considered are $E  =20 $ and $\nu = 0.3$.  Symmetric (traction) boundary conditions are applied as illustrated in Fig.~\ref{Fig_aircraft_mesh}; the remaining boundary segments are subjected to homogenous Dirichlet boundary conditions.  The problem was solved using i-TFEM and FEM. While FEM required 1.85 seconds  to solve the problem, i-TFEM required 2.13  seconds.

To compare the accuracy of FEM and i-TFEM, the reference solution was obtained by solving the same problem over a very fine quadrilateral mesh with nearly 10,000 elements. The problem was then solved over the tangled quad mesh shown in Fig.~\ref{Fig_aircraft_mesh} using FEM and i-TFEM. The i-TFEM solution at the re-entrant vertex had a relative error of 0.018\%, while the error for FEM was 0.11\%. The i-TFEM post-processed solution is illustrated in Fig.~\ref{Fig_aircraft_pp}.

\section{Conclusions}

It is well-known that concave elements lead to a non-invertible parametric mapping which, in turn, leads to erroneous solutions  in standard FEM. To resolve this, an isoparametric tangled finite element method (i-TFEM) was proposed for nonlinear elasticity with material nonlinearity, over bi-linear quadrilateral meshes. 

The proposed method replaces the full-invertibility requirement with the partial invertibility by incorporating compatibility constraints, thereby allowing for tangled (concave) elements to be present in the mesh. The constraints are imposed via Lagrange multipliers, and the corresponding mixed variational formulation was presented. Moreover, the proposed definition of the field in the tangled region naturally leads to modification in the elemental stiffness matrices over  the concave elements. 

Numerical simulations reproducing a set of well-known benchmark tests are promising. i-TFEM converges to the expected solution even in the presence of severely tangled elements. Moreover, the convergence rate and condition number of i-TFEM are comparable over tangled meshes with that of the standard FEM over non-tangled meshes. 

The proposed method can potentially be extended to three dimensions and isogeometric analysis where tangling is known to occur \cite{fusseder2015fundamental, xia2018generating}. 

%i-TFEM is based on the fundamental idea of separating the piecewise invertible mapping. 

%Thus, the method eases the convexity restriction on the FE mesh generators and entails minimal modifications in the exisiting FE framework.

%In this paper, i-TFEM is presented to handle large deformation problems with material non-linearity. The proposed method requires minimal changes to the standard nonlinear FEM framework, and numerical experiments demonstrate its efficiency and convergence characteristics. The results showed that i-TFEM can handle tangled meshes with the same convergence rate as the standard FEM with regular meshes, with the condition number of the stiffness matrix growing at the same rate as standard FEM with untangled/regular mesh. 

\section*{Acknowledgments}
	The authors would like to thank the support of National Science Foundation through grant 1715970.
\bibliography{mybibfile, TFEM3D_new}

\end{document}